\documentstyle[12pt,psfig]{article}

\def\figura#1#2{#2} 	

\reversemarginpar

\makeatletter
\def\@begintheorem#1#2{\it \trivlist \item[\hskip \labelsep{\bf #1\ #2.}]}
\newtheorem{teo}{Theorem}[section]
\newtheorem{rem}[teo]{Remark}
\newtheorem{lem}[teo]{Lemma}

\newtheorem{cor}[teo]{Corollary}

\newtheorem{prop}[teo]{Proposition} 
\newtheorem{que}[teo]{Question} 
\def\finedim#1{{\hfill\hbox{\enspace\fbox{\ref{#1}}}}\vspace{5pt}}
\def\dim#1{\vspace{1pt}\noindent{\it Proof of} {\hspace{2pt}}\ref{#1}.}

\def\compo{\,{\scriptstyle\circ}\,}
\def\cont{{\rm C}}
\def\trasvint{\cap\hspace{-8.5pt}|\hspace{5pt}}
\def\interior{{\rm Int}}

\font\scpicc=cmcsc10
\font\titsc=cmcsc10 scaled 1400
\font\sc=cmcsc10 scaled 1200
\newfont{\Bbb}{msbm10 scaled 1200}
\newfont{\mycal}{eusm10 scaled 1200}
\newfont{\Got}{eufm10 scaled 1200}
\newfont{\Gotpicc}{eufm10}
\def\mr{{\hbox{\Bbb R}}}

\def\mz{{\hbox{\Bbb Z}}}

\def\bb{{\cal B}}

\def\ee{{\cal E}}
\def\rr{{\cal R}}

\def\pp{{\cal P}}

\def\dd{{\cal D}}
\def\call{{\cal L}}
\def\lgot{{\hbox{\Got l}}}
\def\pgot{{\hbox{\Got p}}}
\def\fgot{{\hbox{\Got f}}}
\def\lgotpicc{{\hbox{\Gotpicc l}}}
\def\pgotpicc{{\hbox{\Gotpicc p}}}
\def\fgotpicc{{\hbox{\Gotpicc f}}}
\def\hatM{{\widehat{M}}}
\def\hatv{{\widehat{v}}}

\def\ristr#1{\big|_{#1}}

\def\Comb{{\rm Comb}}
\def\Fram{{\rm Fram}}
\def\Leg{{\rm Leg}}
\def\Pleg{{\rm PLeg}}
\def\foi{finite-order invariant}
\def\stwotriv{S^2_{\rm triv}}
\def\bthtr{B^3_{\rm triv}}
\def\hatpp{\widehat{\pp}}

\title{Combed 3-Manifolds with Concave Boundary,
Framed Links, and Pseudo-Legendrian Links}

\author{Riccardo {\titsc Benedetti}\qquad Carlo {\titsc Petronio}}

\begin{document}

\maketitle

\noindent{\small{\scpicc Abstract}. We provide combinatorial realizations,
according to the usual objects/moves scheme, of the following three topological categories:
(1) pairs $(M,v)$ where $M$ is a 3-manifold (up to diffeomorphism)
and $v$ is a (non-singular vector) field, up to homotopy; here possibly $\partial M\neq\emptyset$, and $v$
may be tangent to $\partial M$, but only in a concave fashion, and homotopy should
preserve tangency type; (2) framed links $L$ in $M$, up to
framed isotopy; (3) triples $(M,v,L)$, with $(M,v)$ as above and $L$
transversal to $v$, up to pseudo-Legendrian isotopy
(transversality-preserving
simultaneous homotopy of $v$ and isotopy of $L$). All realizations 
are based on the notion of branched standard spine, and build on 
results previously obtained. Links are encoded by means of diagrams on branched
spines, where the diagram is ${\rm C}^1$ with respect to the branching.
Several motivations for being interested in combinatorial realizations of the
topological categories considered in this paper are given in the introduction.
The encoding of links is suitable for the comparison of the framed and the
pseudo-Legendrian categories, and some applications are given in connection 
with contact structures, torsion and finite-order invariants.
An estension of Trace's notion of winding number of a knot diagram is 
introduced and discussed.}

\vspace{.3cm}


\vspace{.2cm}

\section*{Introduction}
This paper describes combinatorial realizations,
based on the machinery of branched standard spines (see Section~\ref{reminder})
of the following three topological categories (in which manifolds and
diffeomorphisms are oriented by default):

\begin{enumerate}

\item {\em Combed $3$-manifolds with concave boundary}, that is pairs 
$(M,v)$, where $M$ is a compact 3-manifold (possibly with boundary), and $v$ is a nowhere-zero
vector field on $M$ with simple tangency circles of concave type on $\partial M$, up to
diffeomorphism of $M$ and homotopy of $v$ through fields of the same sort;

\item {\em Framed links in $3$-manifolds}, that is 
pairs $(M,L)$, where $M$ is as above and $L$ is a framed link in $M$, 
up to diffeomorphism of $M$ and framed isotopy of $L$;

\item {\em Pseudo-Legendrian links in combed $3$-manifolds}, that is
triples $(M,v,L)$, where $(M,v)$ is as above and $L$ is transversal to $v$,
up to diffeomorphism of $M$ and `pseudo-Legendrian isotopy' of $(v,L)$, 
{\em i.e.}~simultaneous homotopy of $v$ and isotopy of $L$
through pairs $(v,L)$ of the same type.

\end{enumerate}

\noindent  We will denote these categories respectively by $\Comb$,
$\Fram$ and $\Pleg$. (Regarding names, recall that a non-zero vector field up
to homotopy is often called a {\em combing}, and that if $\xi$ is an oriented contact
structure and $L$ is {\em Legendrian} in $\xi$, then $(M,\xi^\perp,L)$ defines
an element of $\Pleg$.) Our realizations are given according to the  by now
popular scheme in 3-dimensional topology, namely:

\begin{enumerate}

\item[(I)] A class of combinatorial objects, each of which can be specified
by a finite set of data, and a surjective 
reconstruction map which assigns to a combinatorial object a
topological one;

\item[(II)] A finite set of local combinatorial moves on objects,
finite combinations of which give the equivalence relation induced by the
reconstruction map.

\end{enumerate}

In the definitions of the topological categories given above we have been
forced to include the action of diffeomorphisms, because we use spines, which
determine manifolds only up to diffeomorphism. However if a certain manifold
$M$ is given we can restrict to spines embedded in $M$ (rather than abstract
ones), and get formally identical combinatorial realizations of the refined
categories where only diffeomorphisms of $M$ isotopic to the identity are
considered. We will mention how to do this in Section~\ref{conc:calc} for
combings, but a similar refinement could easily be stated for framed
links and for pseudo-Legendrian links.

Rather than providing precise statements of our realizations, in this introduction
we give some general background and motivations, starting with $\Comb$.
A combinatorial realization of the subcategory $\Comb^{\rm cl}$
of $\Comb$ given by pairs $(M,v)$ with 
closed $M$ was given in~\cite{lnm}. In Section~\ref{conc:calc} we extend the arguments
of~\cite{lnm} to the bounded case, and we actually refine the results proved there, by
showing that some of the moves previously considered may actually be neglected.
The realization of $\Comb^{\rm cl}$ in \cite{lnm} was the basis for the 
treatment of other refinements of the category of 3-manifolds, involving
spin structures and framings. These realizations proved fruitful in connection
with spin-refined Turaev-Viro invariants (see Section~8.3
in~\cite{lnm}) and G.~Kuperberg's invariants for combed and framed manifolds,
of which a very constructive description is given in~\cite{benvenuti}.
Our main motivation here
comes from~\cite{first:paper}, where we have developed a theory of
Euler structures with simple boundary and their Reidemeister-Turaev torsion
(see~\cite{turaev:Euler},~\cite{turaev:spinc},~\cite{turaev:nuovo}).
The surjectivity of the reconstruction map of the realization of $\Comb$
was used in~\cite{first:paper}
to construct an explicit canonical $H_1$-equivariant bijection
from the space of smooth Euler structures to the space of
combinatorial Euler structures, and to exhibit a canonical Euler chain
for the structure carried by a branched spine.

The subcategory of $\Fram$ consisting of framed links in closed manifolds
was combinatorially realized by Turaev~\cite{turaev:ombre}
in terms of link diagrams on a given standard spine, and 
moves on these diagrams (including the classical framed Reidemeister moves). 
In Section~\ref{fram:calc} we modify the situation considered by Turaev by
taking a {\em branched} standard spine of the manifold, and restricting to
link diagrams which are $\cont^1$ with respect to the branching. On one hand, 
this allows to simplify
the encoding of the framing, because the field carried by the spine is automatically
transverse to the link, while Turaev needs to add half-twists. On the other hand,
some technical complications emerge, 
because only $\cont^1$ moves can be used. Nevertheless,
a result formally analogous to Turaev's turns out to be true, yielding
the presentation of $\Fram$ discussed in Section~\ref{fram:calc}.
In Section~\ref{leg:calc} we exploit the fact that if a 
branched spine defines a global field on a manifold, according to the
scheme given for $\Comb$, then
the link defined by a $\cont^1$ diagram on the spine 
is automatically pseudo-Legendrian with respect to the field. 
This leads us to the presentation of $\Pleg$.

Comparing the presentations of $\Fram$ and $\Pleg$ one notices a
rather remarkable feature: the former is obtained from the latter just by adding
the `curl' (first Reidemeister) move. This fact has two interesting interpretations:
\begin{itemize}
\item[(a)] it is a perfect combinatorial
analogue of the imitation of a framed isotopy by a Legendrian isotopy in a 
contact manifold;
\item[(b)] it allows a partial extension of the notion 
of winding number of a link diagram.
\end{itemize}
The imitation mentioned in (a)
plays a central role in the comparison, due to Fuchs-Tabachnikov~\cite{futa}
and Tchernov~\cite{vlad} of framed and Legendrian \foi s,
and we believe that our combinatorial realizations could be of some help in the
understanding of these invariants. In particular, we conjecture that the right environment
in which \foi s should be considered in precisely our category $\Pleg$ (we will
provide an exact statement and some evidence in Section~\ref{specu:section}).

Concerning (b), recall first~\cite{trace} that
if a knot $K:S^1\to\mr^3$ is transverse to the constant vertical field $\partial/\partial z$,
then its equivalence class up to isotopy transverse to $\partial/\partial z$ is
determined by the framed isotopy class and by the `winding' number 
(the degree of
$\pi\compo K'$, where $\pi$ is the obvious projection on the horizontal unit circle).
Using our presentations of $\Fram$ and $\Comb$ we can show that
a partial analogue of this fact is true in any combed manifold with concave
boundary, provided one allows a homotopy of the field simultaneous with the isotopy of
the knot. In the general setting, however, the winding number only
exists as a relative object, and we can prove that it leads to a well-defined
invariant only under the assumption that the knot is `good.' The notion of `goodness'
for knots emerged in our study of torsion 
as a relative invariant of pairs of pseudo-Legendrian
knots which are framed-isotopic~\cite{first:paper}. Many knots are good: for instance,
all knots are good if the ambient manifold is a homology sphere, and
most knots with hyperbolic complement are good. 
In Section~\ref{specu:section} we will give some applications of the notion
of relative winding number, in connection with torsion and finite-order invariants.
In particular we will show the following:
\begin{prop}\label{appl:wind:prop}
Let $M$ be a homology sphere, let $v$ be a field on $M$ and consider
two pseudo-Legendrian knots in $(M,v)$ which are
isotopic as framed knots. Then the following conditions are pairwise equivalent:
\begin{enumerate}
\item\label{iso:intro:point} the knots are pseudo-Legendrian isotopic;
\item\label{zero:wind:intro:point} the relative winding number vanishes;
\item\label{same:maslov:intro:point} the knots have the same Maslov index;
\item\label{same:torsion:intro:point} the knots cannot be distinguished by the relative
torsion invariants of~\cite{first:paper};
\item\label{imm:intro:point} the knots are homotopic as pseudo-Legendrian {\em immersions}.
\end{enumerate}
\end{prop}
Moreover we will prove that torsion invariants cannot distinguish the pairs
of framed-isotopic Legendrian knots 
given in~\cite{vlad}, which Tchernov shows to be distinguished by \foi s.

We conclude by giving another perspective of the realization of $\Pleg$.
Recall that $\Pleg$ comes as a refinement of 
Turaev's presentation of $\Fram$, which was the starting point
of his beautiful theory of 4-dimensional shadows. We believe that the extra structure
given by the branching of the spine, which underlies the presentation of $\Pleg$,
should have 4-dimensional counterparts. Our intuition is that
``4-dimensional branched shadows'', which are not quite defined yet,
should correspond to ${\rm Spin}^{\rm c}$ structures on 4-manifolds and
allow to treat their invariants. This intuition is supported by the
fact that in dimension three branched spines indeed are a good framework
to treat torsion of Euler structures (and hence, in particular, Seiberg-Witten
invariants of closed 3-manifolds with ${\rm Spin}^{\rm c}$ structures,
see~\cite{turaev:spinc},~\cite{turaev:nuovo},~\cite{meta}).

\medskip

{\sc Acknowledgment}.
Section~\ref{specu:section} owes a lot to very useful discussions we had with Vladimir
Tchernov.

\section{Branched spines and combings}\label{reminder}
This section contains many definitions used below and
reviews the theory developed in~\cite{lnm}. 

\begin{figure}
\figura{\vspace{.7in}
\centerline{\special{wmf:SecondFigs/SecondWMF/cvtang.wmf y=.7in}\hspace{5cm}\ }}
{\centerline{\psfig{file=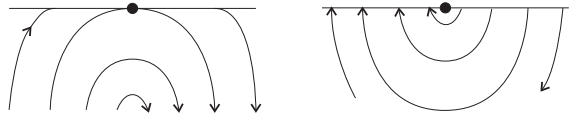,width=7.5cm}}}
\caption{\label{conc:conv:tang}Convex (left) and concave (right) tangency to the boundary.}
\end{figure}

\paragraph{Manifolds and fields}
All the manifolds we will consider are 3-dimensional, oriented, and compact, with
or without boundary. Using the {\it Hauptvermutung}, we will somewhat intermingle the
differentiable and piecewise linear viewpoints. Maps
will always respect orientations. All vector fields mentioned in this paper will be 
non-singular, and they will be termed just {\em fields} for the sake of brevity.
A field $v$ on a manifold $M$ is called {\em traversing} 
if its orbits eventually intersect $\partial M$ 
transversely in both directions (in other words, orbits are compact intervals or points).
A point where $v$ is tangent to $\partial M$ is called {\em simple}
if it appears in a cross-section as in Fig.~\ref{conc:conv:tang}. 
The field is called {\em concave} if it is tangent to $\partial M$ only in a 
concave fashion, as shown on the left in the figure.
Given a concave field $v$ on $M$, the boundary of $M$ naturally splits into 
the region on which $v$ points outside $M$ (which we denote by $B$ and call the 
{\em black} region), and the region on which $v$ points inside
(denoted by $W$ and called {\em white}).
Note that $\partial B=\partial W$ is a union of circles. The pair $(B,W)$,
which is actually determined by any two of its elements, is called a {\em boundary pattern}
on $M$. (This definition is a simplified version of that given in~\cite{first:paper},
because here we do not allow convex tangency.) Starting from the Poincar\'e-Hopf
theorem one can show that a given boundary pattern $\pp=(B,W)$ on $M$, 
{\em i.e.}~a splitting of $\partial M$ into two surfaces with common boundary,
arises from a concave field if and only if $\chi(W)=\chi(M)$. See~\cite{first:paper}.

\paragraph{Standard spines}
A {\em simple} polyhedron $P$ is a  finite, connected, purely 2-dimensional polyhedron
with singularity of stable nature (triple lines and points where six non-singular
components meet). Such a $P$ is called {\it standard} if all the components of
the natural stratification given by singularity are open cells. Depending on
dimension, we will call the components {\it vertices, edges} and {\it regions}.

\begin{figure}
\figura{\vspace{1.5cm}
\centerline{\special{wmf:SecondFigs/SecondWMF/screwbra.wmf y=1in}\hspace{12cm}\ }}
{\centerline{\psfig{file=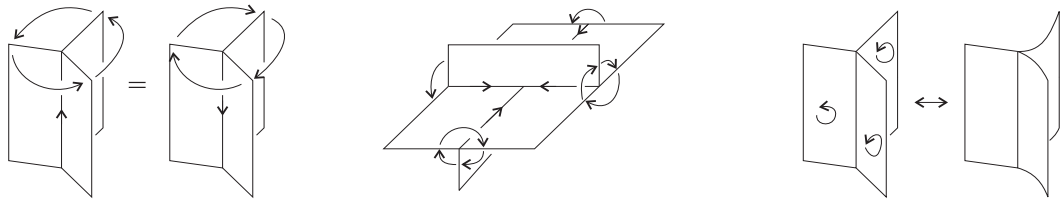,width=14cm}}}
\caption{\label{screw:branch}Convention on screw-orientations, compatibility 
at vertices, and geometric interpretation of branching.}
\end{figure}

A {\em standard spine} of a $3$-manifold $M$ with $\partial M\neq\emptyset$
is a standard polyhedron $P$ embedded in $\interior(M)$ so that $M$ collapses onto $P$.
Standard spines of oriented $3$-manifolds are characterized among standard polyhedra
by the property of carrying an {\em orientation}, defined 
(see Definition~2.1.1 in~\cite{lnm}) as a ``screw-orientation''
along the edges (as in the left-hand-side of Fig.~\ref{screw:branch}),
with the obvious compatibility at vertices
(as in the centre of Fig.~\ref{screw:branch}).
It is the starting point of the theory of standard spines
that every oriented $3$-manifold $M$ with $\partial M\neq\emptyset$ has an oriented
standard spine, and can be reconstructed (uniquely up to equivalence) 
from any of its oriented standard spines. See~\cite{casler} for the non-oriented version 
of this result and~\cite{manuscripta} or Proposition~2.1.2 in~\cite{lnm} for the (slight) oriented refinement. We will denote by $M(P)$ the manifold defined by $P$.
Note that $\partial M(P)\neq\emptyset$. To recover closed manifolds one considers spines $P$
such that $\partial M(P)\cong S^2$, and defines $\hatM(P)$ as $M(P)\sqcup_f D^3$
with $f:\partial M(P)\to S^2$ a diffeomorphism. Note that this definition makes sense
also when $\partial M(P)$ has more than one component, but at least one is a sphere.

\begin{figure}
\figura{\vspace{3cm}}
{\centerline{\psfig{file=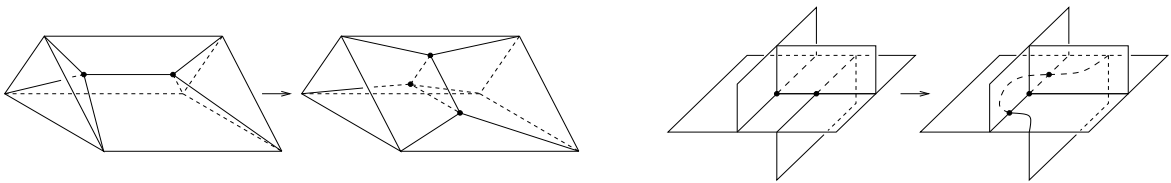,width=15.5cm}}}
\caption{\label{mapimove}Two pictures of the Matveev-Piergallini move.}
\end{figure}

\begin{figure}
\figura{\vspace{3cm}}
{\centerline{\psfig{file=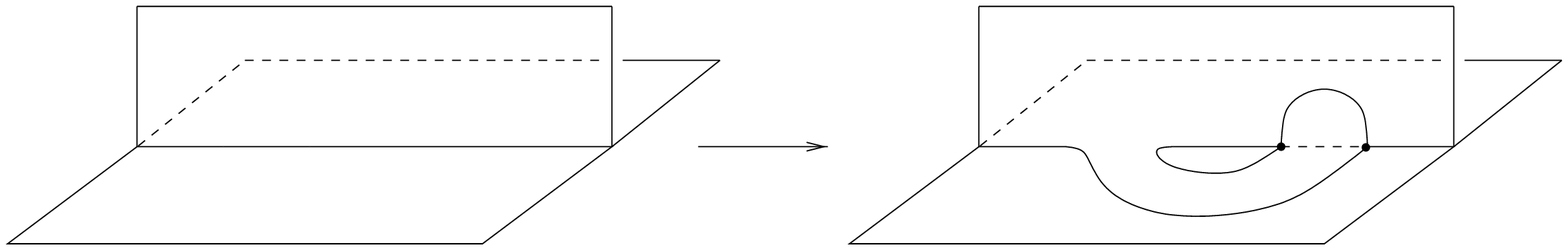,width=10cm}}}
\caption{\label{02move}Another move on standard spines.}
\end{figure}

\paragraph{Moves for standard spines}
The fundamental move for standard spines, which (in both directions) preserves 
the topological
type of the associated manifold, is the Matveev-Piergallini MP
move, see~\cite{matv:mossa},~\cite{piergallini} and Fig.~\ref{mapimove}. 
Counting the vertices involved one is naturally led to 
call the positive MP a ``2-to-3'' move.
The MP-move and its inverse are actually not sufficient to relate spines of the same
manifold, because they obviously cannot apply to spines with one vertex. However,
as soon as one decides to dismiss these ``MP-rigid'' spines (not the corresponding
manifolds, which have plenty of other spines), the MP-move does become 
sufficient~\cite{piergallini}. To deal with spines with one vertex the ``0-to-2''
move of Fig.~\ref{02move} (and its inverse)
must be added.

\paragraph{Branched spines}
A {\it branching} on a standard polyhedron $P$ is an
orientation for each region of $P$, such that no edge is induced the same
orientation three times. See the right-hand side of Fig.~\ref{screw:branch}
and Definition 3.1.1 in~\cite{lnm} for the geometric meaning of this notion.
An oriented standard spine $P$ endowed with a branching is shortly named 
{\em branched spine}. We will never use specific notations for the extra structures: 
they will be considered to be part of $P$.
The following result, proved as Theorem~4.1.9 in~\cite{lnm}, is the starting point of our constructions.

\begin{prop}\label{from:spine:to:field}
To every branched spine $P$ there corresponds a manifold $M(P)$ 
with non-empty boundary and a concave traversing field $v(P)$ on $M(P)$.
The pair $(M(P),v(P))$ is well-defined up to equivalence, and
an embedding $i:P\to\interior(M(P))$ is defined
with the property that $v(P)$ is positively transversal to $i(P)$.
\end{prop}

\begin{figure}
\figura{\vspace{3cm}}
{\centerline{\psfig{file=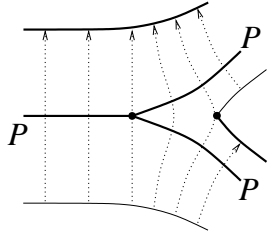,width=3.5cm}}}
\caption{\label{constr:M}Manifold and field associated to a branched spine.}
\end{figure}

The topological construction which underlies this proposition is actually quite
simple, and it is illustrated in Fig.~\ref{constr:M}. Concerning the last
assertion of the proposition, note that the branching allows to define an
oriented tangent plane at each point of $P$.

\paragraph{Non-traversing fields and closed manifolds}
As noted above, standard spines do not directly represent closed manifolds,
but one can use spines of manifolds bounded by $S^2$ and cap off this sphere
to get a closed manifold, or, viewing things the other way around, one
can remove an open ball from a given closed manifold to get a bounded one.
When one is interested in a manifold equipped with a field,
one can try to use branched spines, but of course one sees that they are inadequate 
to give a direct description both when the manifold is closed and when the field
is non-traversing. This limitation is circumvented again by removing a ball, with
a {\em proviso} on the field on that ball.

Let $P$ be a branched standard spine, and assume that in 
$\partial M(P)$ there is only one component which is diffeomorphic to $S^2$
and is split by the tangency line of $v(P)$ to $\partial M(P)$ into two
discs. Such a component will be denoted by $\stwotriv$. Now, notice that $\stwotriv$ is
also the boundary of the closed $3$-ball with constant vertical field, 
denoted by $B^3_{\rm triv}$. This shows that we can cap off
$\stwotriv$ by attaching a copy of $B^3_{\rm triv}$, 
getting a compact manifold $\hatM(P)$ and a 
concave field $\hatv(P)$ on $\hatM(P)$. 
If we denote by $\hatpp(P)$ the boundary pattern of $\hatv(P)$ on
$\hatM(P)$, we easily see that the pair $(\hatM(P),\hatv(P))$ is only well-defined
up to diffeomorphism of $\hatM(P)$ and homotopy of $\hatv(P)$ through fields compatible
with $\hatpp(P)$.

\begin{figure}
\figura
{\vspace{3cm}}
{\centerline{\psfig{file=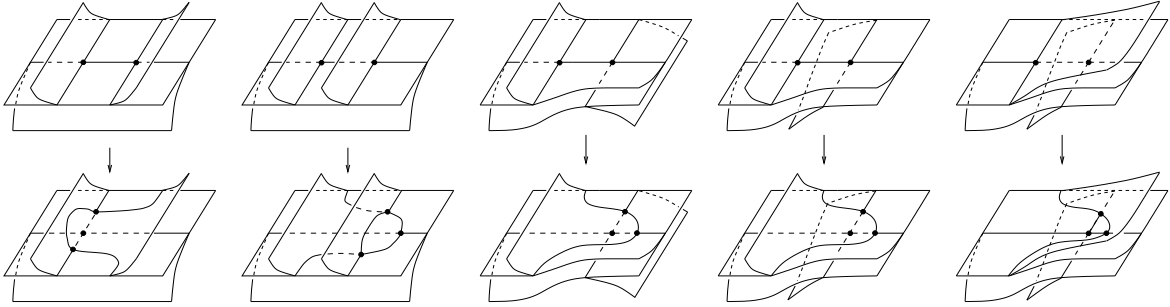,width=15.5cm}}}
\caption{\label{slidingMP}Sliding-MP-moves.}
\end{figure}

\begin{figure}
\figura{\vspace{3cm}}
{\centerline{\psfig{file=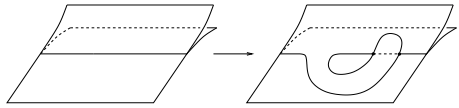,width=6cm}}}
\caption{\label{sliding02}Another sliding move.}
\end{figure}

\paragraph{Standard sliding moves} 
Let $P\to P'$ be a {\em positive} MP-move (so, $P'$ has one
vertex and one region more than $P$). If $P$ has a branching, all the regions of $P'$, 
except for the new one, already have an orientation, and it is a fact that the new
region can always be given an orientation (sometimes not a unique one) so to get a
branching on $P'$. Each of the moves on branched spines arising like this will
be called a {\em branched} MP-move, and it
will be called a {\em sliding} MP-move if moreover it does not modify the
boundary pattern of the associated concave field. One can actually see that
each sliding-MP-move can be realized within a certain pair $(M,v)$ as a
continuous deformation through branched spines of $M$ transverse to $v$, 
with only one singularity at which the spine is non-standard but transversality 
is preserved. This deformation is shown in Fig.~\ref{slidingMP},
and it justifies the term `sliding' quite clearly.
Since in Fig.~\ref{slidingMP} we are showing portions of spines embedded in $\mr^3$,
to give a completely intrinsic description of the moves we should specify in each
portion whether the screw-orientation of the spine is equal or opposite to
that induced by $\mr^3$, and whether the upward vertical field is positively
or negatively transversal to the spine. As a result, the complete list of 
sliding-MP-moves contains 16 different ones, but the essential physical
modifications are only those shown in Fig.~\ref{slidingMP}.
From this figure one also sees quite clearly that if $\hatv(P)$ and $\hatv(P')$
can be defined then they coincide (up to homotopy through fields compatible
with $\hatpp(P)=\hatpp(P')$). Another move which obviously has the same
property, and will be needed below, is the branched version of
the 0-to-2 move, shown Fig.~\ref{sliding02}, and called the {\em snake move}
in the sequel. As above, if one takes orientations into account, 
there is another essentially different snake move, obtained by mirroring
Fig.~\ref{sliding02}. Since also the snake move involves a sliding, we will
call {\em standard sliding move} any sliding-MP or snake move.

\section{A calculus for combed manifolds\\ with concave boundary}\label{conc:calc}
In this section we will extend and refine the main results of Chapter~5 of~\cite{lnm}.
The extension consists in passing from the closed to the bounded case, and the 
refinement comes from the shortening of the list of moves to be considered.
More precisely, we will show that compact manifolds
with concave combings are
combinatorially described by (suitable) branched spines up to certain moves,
namely the standard sliding
(snake and sliding-MP) moves shown above. 
Moreover, we will show that spines which
are rigid with respect to sliding-MP-moves can be dismissed with no harm, and that
the sliding-MP-moves suffice to generate the equivalence on the remaining spines.
This implies that our result is a perfect combed analogue of the Matveev-Piergallini
theorem (but our proof is self-contained).

\paragraph{Definitions and statements}
We will denote by $\Comb$ the set of all
pairs $(M,v)$, where $M$ is a compact oriented manifold and 
$v$ is a concave field on $M$,
viewed up to diffeomorphism of $M$ and homotopy of $v$
through concave fields. A class
$[M,v]\in\Comb$ is called a {\it combing} on the diffeomorphism class of the
manifold $M$. Note that the boundary pattern on $\partial M$
evolves isotopically during a homotopy of $v$, so a pair $(M,\pp)$, viewed
up to diffeomorphism of $M$, can be
associated to each $[M,v]\in\Comb$. In particular, $\Comb$ naturally splits
as the disjoint union of subsets $\Comb([M,\pp])$, consisting of combings on $M$ compatible
with $\pp$.

For a technical reason we actually rule out from $\Comb$ the set of those
classes $[M,v]$ such that the corresponding boundary pattern contains components of the type
$\stwotriv$. This is actually not a serious restriction, because each
$\stwotriv$ component can be capped off by a $\bthtr$, 
and the result is well-defined up to homotopy. Note that we do accept
pairs $(M,v)$ with closed $M$, and pairs in which $v$ has no tangency at all to
$\partial M$.

Let us denote now by $\bb$ the set of all branched spines $P$ (up to PL isomorphism) 
such that $\pp(P)$ contains only one $\stwotriv$. Such a $P$ being given,
$\hatM(P)$ and $\hatv(P)$ can be considered, and the pair
$(\hatM(P),\hatv(P))$ gives rise to a well-defined element of $\Comb$,
which we denote by $C(P)$. The following will be shown below:

\begin{teo}\label{comb:calc:teo}
The map $C:\bb\to\Comb$ is surjective, and the equivalence
relation defined by $C$ on $\bb$ is generated by sliding-MP-moves
and snake moves.
\end{teo}

\begin{rem}\label{trivial:pieces}
{\em The following interpretation of the surjectivity of $C$ is perhaps useful.
Note first that the dynamics of a field, even a concave one, can be very complicated,
whereas the dynamics of a traversing field (in particular, $B^3_{\rm triv}$) is simple.
Surjectivity of $C$ means that for any (complicated) concave field
there exists a sphere $S^2$ which splits the field into two (simple) pieces:
a standard $B^3_{\rm triv}$ and a concave traversing field. Actually, a 1-parameter version
of this statement also holds (see Remark~\ref{one-parameter}): we will need it to show that the 
$C$-equivalence is the same as the sliding equivalence.}
\end{rem}

As announced, we state now the sliding analogue of the fact that the MP moves
suffice. Let us denote by $\rr$ the subset of $\bb$ consisting of the branched spines 
which are ``rigid'' from the point of view of the sliding-{\rm MP}-moves, {\it
i.e.}~the spines to which no such move applies. An explicit description of $\rr$
is given in the proof of the next result. In the statement we only emphasize the
most important consequences of this description. 

\begin{prop}\label{rigid:spines}
\begin{enumerate}
\item[(i)] For every surface $\Sigma$ and pattern $\pp$ on $\Sigma$
there are at most two spines
$P\in\rr$ such that $\partial(M(P))\cong(\Sigma,\pp)$.
\item[(ii)] If two elements of $\bb\setminus\rr$ are related through
sliding-{\rm MP}moves and snake moves, they are also related through 
sliding-{\rm MP}-moves only.
\item[(iii)] Every $P\in\rr$ is related by a snake move to
an element of $\bb\setminus\rr$.
\end{enumerate}
\end{prop}

This proposition shows that in the statement of Theorem~\ref{comb:calc:teo}
one may remove $\rr$ from $\bb$ and forget the snake move.

\paragraph{Embedding-refined calculus}
We spell out in this paragraph the embedding-refined version of our calculus,
which allows to neglect the action of automorphisms. Let a certain manifold $M$ be given, and consider the set $\Comb(M)$ of
concave vector fields on $M$, up to homotopy. Let $\bb(M)$ consist of the
elements of $\bb$ which are smoothly embedded in $M$ as spines of $M$ minus a ball
$\bthtr$. Each element $P$ of $\bb(M)$ is viewed up to isotopy in $M$, and gives rise to
a well-defined element $C_M(P)$ of $\Comb(M)$. Moreover sliding-MP-moves
and snake moves are well-defined in $\bb(M)$, because they can be realized as embedded
moves. The embedded analogue of Theorem~\ref{comb:calc:teo} states that
$C_M:\bb(M)\to\Comb(M)$ is surjective, and the relation it defines is generated by the
embedded moves. The proof of this result is a refinement of the 
proof of the general statement, along the lines explained
in~\cite{lnm} (4.1.12, 4.1.13, 4.3.5, and~5.2.1.)

\paragraph{Normal sections of a concave field}
The proof of Theorem~\ref{comb:calc:teo} is an extension of the argument
given in Chapter~5 of~\cite{lnm}, and it is based on the following technical
notion, which extends ideas originally due to Ishii~\cite{ishii}. Let $v$ be a
concave field on $M$. Let $B_1,\dots,B_k$ be the black components of the
splitting of $\partial M$, {\it i.e.}~the regions on which $v$ points outwards. A
{\em normal section} for $(M,v)$ is a compact surface $\Sigma$ with boundary,
embedded in the interior of $M$, with the following properties:
\begin{enumerate}
\item\label{transverse:point} $v$ is transverse to $\Sigma$;
\item\label{disc:point} $\Sigma$ has exactly $k+1$ components $\Sigma_0,\dots,\Sigma_k$, with 
$\Sigma_0\cong D^2$;
\item\label{black:point} For $i>0$, the projection of $B_i$ on $\Sigma$ along the orbits of $-v$ is
well-defined and yields a diffeomorphism between $B_i$ and a surface $B'_i$
contained in the interior of $\Sigma_i$, with $\Sigma_i\setminus B'_i$ being a
collar on $\partial\Sigma_i$ (so $B'_i=\Sigma_i$ if $\partial\Sigma_i=\emptyset$);
\item\label{all:orbits:met} Each positive half-orbit of $v$ meets either the interior of some $B_i$
(where it stops), or the interior of some $\Sigma_i$;
\item\label{generic:point} $\partial\Sigma$ meets itself generically along $v$ 
({\it i.e.}~each orbit of
$v$ meets $\Sigma$ at most two consecutive times on $\partial\Sigma$, and, if
so, transversely);
\item\label{standardness:point} Let $P_\Sigma$ be the union of $\Sigma$ with all the orbit
segments starting on $\partial\Sigma$ and ending on $\Sigma$. Then $\Sigma$,
which is a simple polyhedron by the previous point, is actually
standard.
\end{enumerate}

The next two lemmas show that normal sections of $(M,v)$ correspond bijectively
to spines $P$ such that $C(P)=[M,v]$. The proof of surjectivity of $C$
and the discussion of its non-injectivity will be based on these lemmas.

\begin{lem}\label{section:then:spine}
If $(M,v)$, $\Sigma$ and $P_\Sigma$ are as above, then $P_\Sigma$ can be given
a structure of branched spine such that $C([P_\Sigma])=[M,v]$.
\end{lem}

\dim{section:then:spine}
We orient $\Sigma$ so that $v\trasvint\!^+\,\Sigma$
(by default $M$ is oriented). Every region of $P_\Sigma$
contains some open portion of $\Sigma$, so it can be oriented accordingly;
with the obvious screw-orientation, this turns $P_\Sigma$ into a branched spine
of its regular neighbourhood in $M$.

We show that $C([P_\Sigma])=[M,v]$ by embedding the abstract  
manifold $M(P_\Sigma)$ in $M$, in such a way that the field carried by
$P_\Sigma$ on $M(P_\Sigma)\subset M$ is just the restriction of $v$. By
construction, $M\setminus M(P_\Sigma)$ will
consist of a copy of $\bthtr$, together with a collar on $\partial M$ which can
be parameterized as $(\partial M)\times[0,1]$ in such a way that $v$ is constant
in the $[0,1]$-direction. This easily implies that $C([P_\Sigma])=[M,v]$
indeed.

\begin{figure}
\figura{\vspace{3cm}}{
\centerline{\psfig{file=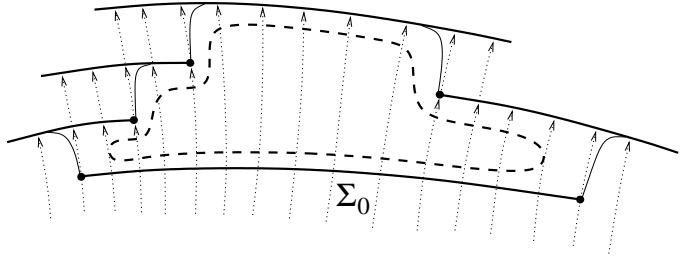,width=9cm}}}
\caption{\label{trivball}The trivial ball.}
\end{figure}

We illustrate the embedding of $M(P_\Sigma)$ in $M$ 
pictorially in one dimension less.
Figure~\ref{trivball} shows how $\Sigma_0$ gives rise to a $\bthtr$. 
In the figure we describe $v$ by dotted lines,
$\Sigma$ by thick lines, portions of $P_\Sigma\setminus\Sigma$ by thin lines,
and $\partial(M(P_\Sigma))$ by a thick dashed line. Note also that the portions
of $P_\Sigma\setminus\Sigma$ have been slightly modified so to become positively
transversal to $v$, which allows us to represent the branching as usual,
{\it i.e.}~as a $\cont^1$ structure on $P_\Sigma$.

\begin{figure}
\figura{\vspace{3cm}}{
\centerline{\psfig{file=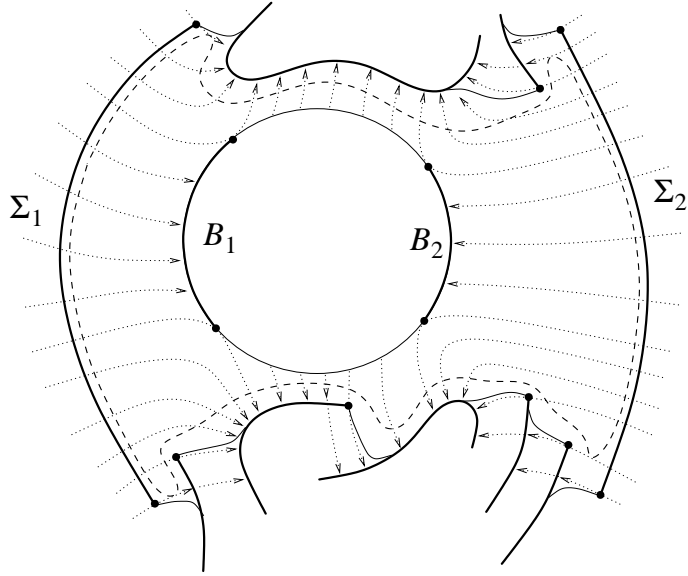,width=9cm}}}
\caption{\label{collar}Collar on a boundary component.}
\end{figure}

Figure~\ref{collar} shows the collar based on a component of $\partial M$.
We use the same conventions as in the previous figure, and in addition we
represent the black and white components of $\partial M$ by thick and thin lines
respectively. This description concludes the proof.
\finedim{section:then:spine}

\begin{lem}\label{spine:then:section}
Let $[P]\in\bb$ and $C([P])=[M,v]\in\Comb$, with $P$ embedded in 
$(M,v)$ according to the geometric description of $C$. Let $\Sigma$ be
obtained from $P$ as suggested (in one dimension less) in Fig.~\ref{cutspine}.
\begin{figure}
\figura{\vspace{3cm}}{
\centerline{\psfig{file=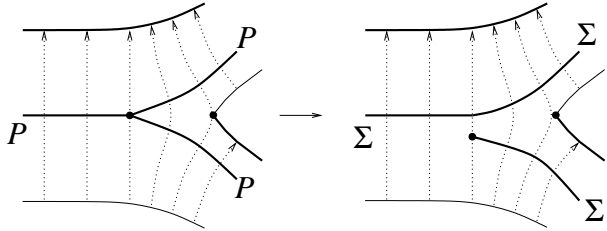,width=8cm}}}
\caption{\label{cutspine}Normal section from a spine.}
\end{figure}
Then $\Sigma$ is a normal section of $(M,v)$, and $P_\Sigma$ is isomorphic to
$\Sigma$.
\end{lem}

\dim{spine:then:section}
The construction suggested by Fig.~\ref{cutspine} is obviously the inverse of 
the construction in the proof of Lemma~\ref{section:then:spine}.\finedim{spine:then:section}

\paragraph{The concave combing calculus} 
Using normal sections we can now
show the main result of this section.

\dim{comb:calc:teo} We start with the proof of surjectivity. So, let
us consider a combed manifold $(M,v)$, subject to the usual restrictions. By
Lemma~\ref{section:then:spine} it is natural to try and construct a normal
section for $(M,v)$. Let $B_1,\dots,B_k$ be the black regions in $\partial M$.
Slightly translate each $B_i$ along $-v$, getting $B'_i$. Add to each $B'_i$ a
small collar normal to $v$, getting $\Sigma_i$ (if $\partial B_i=\emptyset$,
we set $\Sigma_i=B'_i$). Select finitely many discs $\{D_n\}$ disjoint from
each other and from all the $\Sigma_i$'s, such that all positive orbits of $v$,
except for the small segments between $B'_i$ and $B_i$, meet
$(\bigcup_{i\geq1}\Sigma_i)\cup(\bigcup D_n)$ in some interior point. Connect
the $D_n$'s together by strips transversal to $v$ and disjoint from
$\bigcup_{i\geq1}\Sigma_i$, getting a disc $\Sigma_0$. Up to a generic small
perturbation, the surface $\Sigma=\bigcup_{i\geq0}\Sigma_i$ satisfies all
axioms of a normal section for $(M,v)$, except axiom~\ref{standardness:point}.

Now, even if it is not standard, $P_\Sigma$ can be defined, and the proof of
Lemma~\ref{section:then:spine} shows that it is a simple branched spine
of $(M\setminus B^3,v)$. In particular, $P_\Sigma$ is connected and its
singular locus is non-empty. 
We recall now that in Chapter~4 of~\cite{lnm} we have considered a set
of local moves on simple branched spines, called `simple sliding moves',
which preserve the transversal field
(and hence the splitting of the boundary), but do
not require or preserve the cellularity condition.
Knowing that $P_\Sigma$ is connected and $S(P_\Sigma)\neq\emptyset$,
it is not too hard to see
that there exists a sequence of (abstract) simple sliding moves which
turns $P_\Sigma$ into a standard spine (see~\cite{lnm}, Section~4.4). 
If we physically realize these moves
within $M$, preserving transversality to $v$, the result is a standard branched
spine $P$ such that $C([P])=[M,v]$.

We are left to show that if $C([P_0])=C([P_1])$ then $P_0$ and $P_1$
are related by sliding-MP-moves and
snake moves (`sliding-equivalent' for short). 
By the definition of $\Comb$ and $C$, using also
the above lemmas, there exists a manifold $M$ and a homotopy $(v_t)$ of concave
fields on $M$, such that $P_0$ and $P_1$ are defined by normal sections
$\Sigma^{(0)}$ and $\Sigma^{(1)}$ of $(M,v_0)$ and $(M,v_1)$ respectively.

We prove that $P_0$ and $P_1$ are sliding-equivalent first in the special case
where $v_0=v_1=v$. The general case will be an easy consequence. For $j=0,1$,
let $\Sigma^{(j)}=\bigcup_{i\geq0}\Sigma^{(j)}_i$. Proceeding as in the above
proof of surjectivity, for each black region $B_i$ of $\partial M$, we consider
a collared negative translate $\overline\Sigma_i$ of $B_i$. We choose
$\overline\Sigma_i$ so close to $B_i$ that 
$\overline\Sigma_i\cap\Sigma^{(j)}=\emptyset$, and the negative integration of
$v$ yields a diffeomorphism from $\overline\Sigma_i$ to a subset of
$\Sigma^{(j)}_i$. 

\begin{figure}
\figura{\vspace{3cm}}{
\centerline{\psfig{file=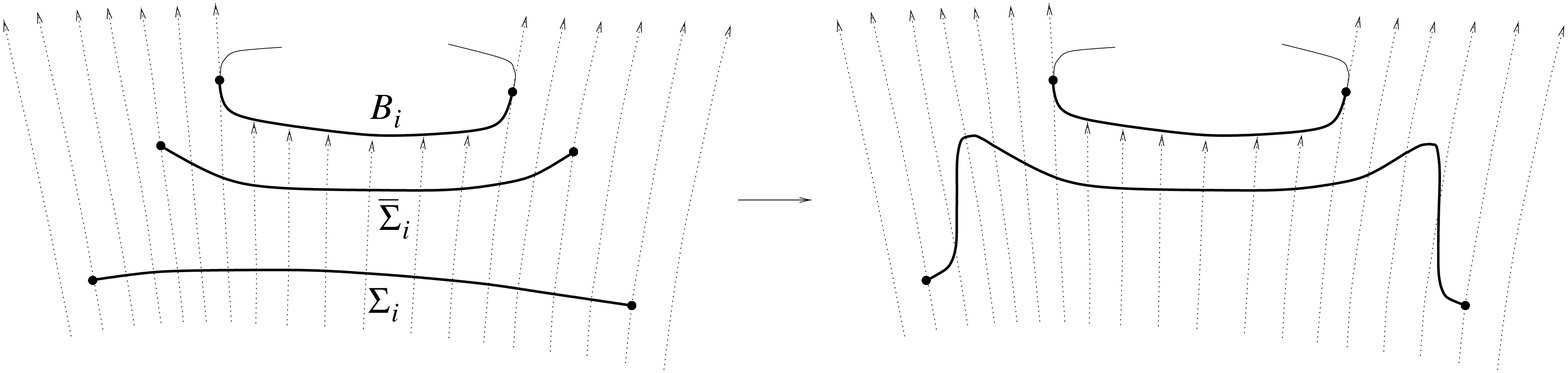,width=15.5cm}}}
\caption{\label{modify:1}Isotopy of a normal section.}
\end{figure}

{\sc Step I}. {\em For $j=0,1$, there exists a disc $D_j$ such that
$D_j\cup(\bigcup_{i\geq1}\overline\Sigma_i)$ is a normal section of $(M,v)$, and
the associated branched spine is sliding-equivalent to $P_j$.} To prove this, we temporarily
drop the index $j$. We first isotope each
$\Sigma_i$, without changing the associated spine, until it contains
$\overline\Sigma_i$, as suggested in Fig.~\ref{modify:1}.

\begin{figure}
\figura{\vspace{3cm}}{
\centerline{\psfig{file=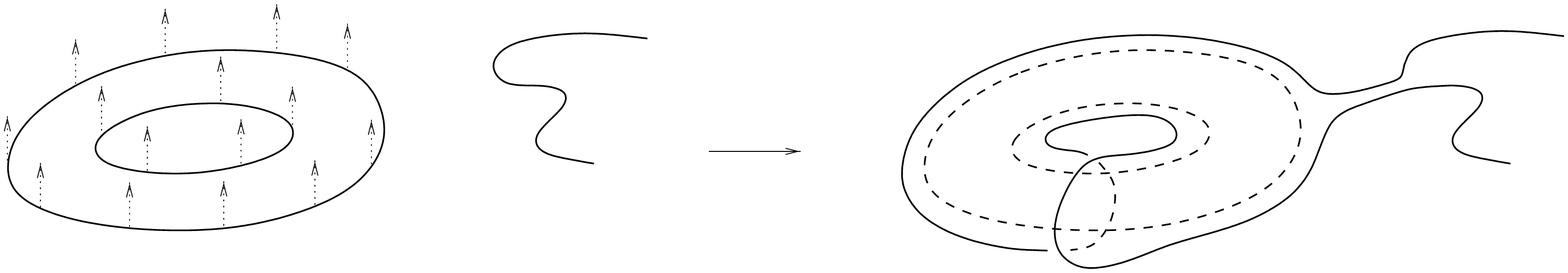,width=14cm}}}
\caption{\label{modify:2}Insulation of an annulus.}
\end{figure}

Note that if $\partial B_i=\emptyset$ we automatically have 
$\Sigma_i=\overline\Sigma_i$. Otherwise, we concentrate on one of the
annuli $A$ of which $\Sigma_i\setminus\overline\Sigma_i$ consists. Note
that we cannot just shrink $A$ leaving the rest of the section
unchanged, because we could spoil axiom~\ref{all:orbits:met} of the definition
of normal section. To actually shrink $A$ we first need to ``insulate''
it, toward the positive direction of $v$, by adding to the disc $\Sigma_0$
a strip normal to $v$. Figure~\ref{modify:2} suggests how to do this.

As we modify $\Sigma_0$ as suggested, it is clear that we keep having a
``quasi-normal'' section, {\em i.e.}~all axioms except~\ref{standardness:point}
hold. Moreover the corresponding simple branched spines are obtained from each
other by the simple sliding moves already mentioned above. 
To conclude we apply, as
above, the fact that a simple branched spine can be transformed via
simple sliding moves to a standard one, and the technical result established
in~\cite{lnm}, Proposition~4.5.6,
according to which standard spines which are equivalent under simple
sliding moves are also sliding-equivalent. This proves
Step I.

The conclusion will now follow quite closely the argument in~\cite{lnm},
Theorem~5.2.1.

{\sc Step II}. {\em There exist discs $D'_0$ and $D'_1$ such that 
$D'_j\cup(\bigcup_{i\geq1}\overline\Sigma_i)$ is a normal section of $(M,v)$
for $j=0,1$, and $D_0\cap D'_0=D'_0\cap D'_1=D'_1\cap D_1=\emptyset$.} Choosing
a metric on $M$, one can construct $D'_0$ and $D'_1$ by first taking many very
small discs almost orthogonal to $v$, and then connecting these discs by strips
transversal to $v$. 

\begin{figure}
\figura{\vspace{3cm}}{
\centerline{\psfig{file=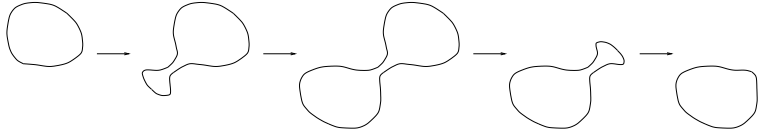,width=10cm}}}
\caption{\label{empty:fill}Transformation of a disc into a disjoint one.}
\end{figure}

{\sc Step III}. {\em Conclusion in the case $v_0=v_1$}. If we connect $D_0$ and
$D'_0$ by a strip orthogonal to $v$, we get a bigger disc $\tilde D_0$ such
that $\tilde D_0\cup(\bigcup_{i\geq1}\overline\Sigma_i)$ is still a normal
section of $(M,v)$. We can actually imagine a dynamical process, in which $D_0$
is first enlarged to $\tilde D_0$, and then is reduced to $D'_0$, as in
Fig.~\ref{empty:fill}.
If the transformation is chosen generic enough, at all times
axioms~\ref{transverse:point},~\ref{disc:point},~\ref{black:point} 
and~\ref{all:orbits:met} will hold, and axiom~\ref{generic:point} will
hold at all but finitely many times. This means that the corresponding branched
spines are related by simple sliding moves. Similarly, we can replace 
$D'_0$ first by $D'_1$ and then by $D_1$.
Using the facts quoted
above, the conclusion follows.

We are left to deal with the general case, where $(v_t)$ is a non-constant
homotopy. It is then sufficient to take a partition $0=t_0<t_1<\dots<t_n=1$ of
$[0,1]$, fine enough that $(M,v_{t_{k-1}})$ and $(M,v_{t_{k}})$ admit a common
normal section which gives rise to isomorphic branched
spines.\finedim{comb:calc:teo}

\begin{rem}\label{one-parameter}
{\em Along the lines of the previous proof we have established
the following topological fact, whose statement does not involve spines. Let
$(v_t)$ be a homotopy of concave fields on $M$, let $B_0,B_1\subset M$ be balls with
$(B_j,v_j)\cong\bthtr$ and $v_j$ traversing on $M\setminus B_j$ for $j=0,1$.
Then there exist another homotopy $(v'_t)$ between $v_0$ and $v_1$ and an
isotopy $(B_t)$ with $(B_t,v_t)\cong\bthtr$ and $v_t$ traversing on $M\setminus
B_t$ for all $t$.}\end{rem}

\paragraph{Sufficiency of the sliding-MP-moves}
To show Proposition~\ref{rigid:spines} we will find it convenient to use the graphic
representation of branched spines introduced in~\cite{lnm}, Section 3.2, but
we do not reproduce here the technicalities needed to introduce this representation.

\begin{figure}
\figura{\vspace{3cm}}{
\centerline{\psfig{file=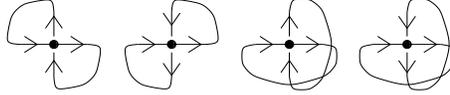,width=6cm}}}
\caption{\label{one:vert:spi}Branched spines with one vertex.}
\end{figure}

\dim{rigid:spines} 
We start by listing rigid spines. Note first that if a negative 
sliding-{\rm MP}-move applies to a spine then also a positive one does, so we only need
to consider positive rigidity. The spines with one vertex, shown in
Fig.~\ref{one:vert:spi}, are of course rigid.
Using~\cite{lnm}, Proposition~3.3.5,
one easily checks that $\partial(M(P))$ is $\stwotriv$ for the
first two spines, and $\stwotriv\sqcup\stwotriv$ for the other two.

\begin{figure}
\figura{\vspace{3cm}}{
\centerline{\psfig{file=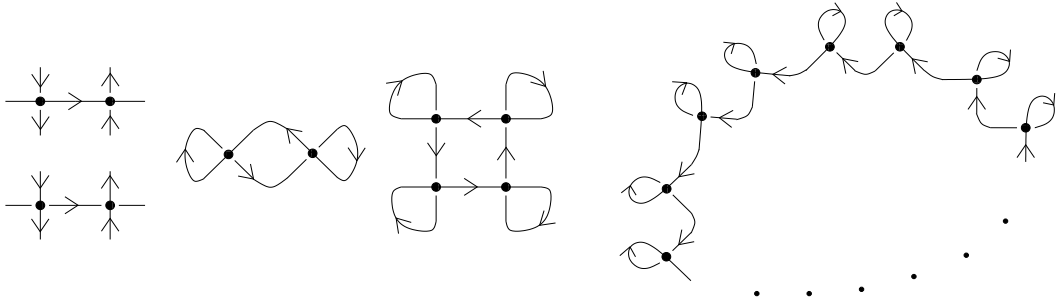,width=14cm}}}
\caption{\label{rigid:spines:fig}Sliding-MP-rigid branched spines.}
\end{figure}

Now we turn to rigid spines with more than one vertex. Rigidity implies that
all edges with distinct endpoints should appear as on the left in 
Fig.~\ref{rigid:spines:fig}.
It is not hard to deduce that rigid spines come in a sequence $P_1^{\rm rig}$,
$P_2^{\rm rig}$, $\dots$ as
shown in the rest of Fig.~\ref{rigid:spines:fig}, where $P_k^{\rm rig}$
has $2k$ vertices, and $\partial(M(P_k^{\rm rig}))$ is the union of
$\stwotriv$ together with $k$ copies of $S^2_{\rm white}$ and 
$k$ copies of $S^2_{\rm black}$.
This classification proves {\em (i)}. 

To show {\em (ii)} we must prove that: 
\begin{enumerate}
\item[{\em (ii-a)}] Sequences which contain rigid spines can be replaced
by sequences which do not.
\item[{\em (ii-b)}] If two non-rigid spines are related by
one snake move then they are also related by a sequence of sliding-{\rm MP}-moves.
\end{enumerate}

For {\em (ii-a)}, we note that the result of a positive snake move is never
rigid. So if a rigid spine $P$ appears in a sequence of moves then $P$ is the
result of a negative snake move
$\mu_1^{-1}:P_1\to P$, and a positive snake move
$\mu_2:P\to P_2$ is applied to $P$.
Since all edges of a spine survive through a snake move,
there is a version $\tilde\mu_2$ of $\mu_2$ which applies to $P_1$ and 
a version $\tilde\mu_1$ of $\mu_1$ which
applies to $P_2$, and the result $\tilde{P}$ is the same. So can
replace the segment $(P_1,P,P_2)$ by $(P_1,\tilde{P},P_2)$, and now all the
spines involved are non-rigid.  

\begin{figure}
\figura{\vspace{3cm}}{
\centerline{\psfig{file=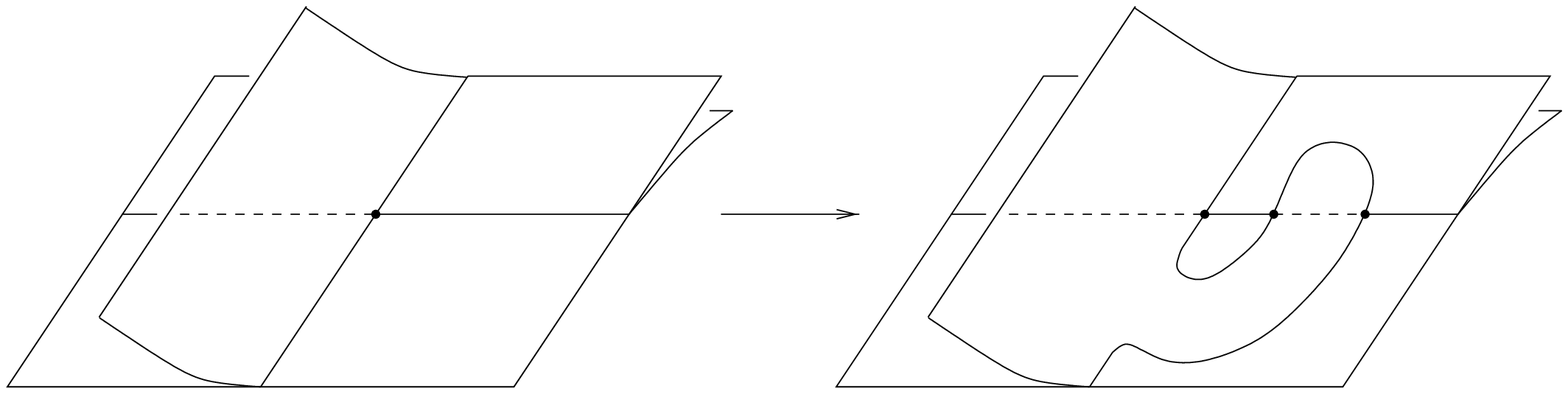,width=7.5cm}}}
\caption{\label{vertex:move}The sliding-vertex move.}
\end{figure}

Let us turn to {\em (ii-b)}. The proof results from three steps, to describe which
we introduce in Figure~\ref{vertex:move}  
another move, called sliding-vertex move, whose unbranched version was already
considered in~\cite{matv:mossa} and~\cite{piergallini}. Again, taking into
account orientations, there are two versions of the move (for each vertex
type), but we will ignore this detail. 

{\em Step 1: if $v$ is a vertex of a branched spine $P$, $e$ is any one of the
edges incident to $v$, $P_v$ is obtained from $P$ via the sliding-vertex move at $v$,
and $P_e$ is obtained from $P$ via the snake move on $e$, then $P_v$ and
$P_e$ are related by sliding-{\rm MP}-moves.} This is proved by an easy case-by-case
analysis. It turns out that two {\rm MP}-moves (a positive and a negative one)
are always sufficient.

{\em Step 2: let $v$, $P$ and $P_v$ be as above. If $P$ and $P_v$  are related
by sliding-{\rm MP}-moves, the same is true for $P$ and any spine obtained from $P$ by
a snake move.} To see this, use step 1 to successively transform sliding-vertex moves
into snake moves and conversely, until the desired snake move is reached.

{\em Step 3: if $P$ is non-rigid then there exists a vertex $v$ such that
$P$ and $P_v$ are related by {\rm MP}-moves.} The vertex $v$ is chosen to be an
endpoint of an edge to which the positive {\rm MP}-move applies. The argument is
again a long case-by-case one, which refines in a branched context the argument
given by Piergallini in~\cite{piergallini}. The sequence always consists of
three positive moves followed by a negative one.
This concludes the proof of {\em (ii)}, whereas {\em (iii)} is evident.\finedim{rigid:spines}

\section{A calculus for framed links}\label{fram:calc}
We fix in this section a compact manifold $M$ and consider the set $\Fram(M)$ of
isotopy classes of framed links in $M$. Since a link isotopy generically avoids a 
fixed 3-ball, $\Fram(M)$ and $\Fram(\hatM)$ are canonically isomorphic
when $\partial M=S^2$, so we can restrict
to non-closed $M$'s and include the closed case as usual. 

\begin{figure}
\figura{\vspace{0.7in}
\centerline{\special{wmf:SecondFigs/SecondWMF/c1proj.wmf y=0.7in}\hspace{6cm}\ }}
{\centerline{\psfig{file=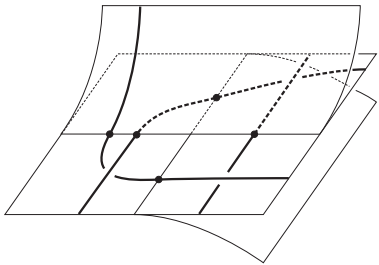,width=5cm}}}
\caption{\label{c1proj}A portion of $\cont^1$ link diagram on a branched spine.}
\end{figure}

\begin{figure}
\figura{\vspace{2.5in}
\centerline{\special{wmf:SecondFigs/SecondWMF/c1moves.wmf y=2.5in}\hspace{10cm}\ }}
{\centerline{\psfig{file=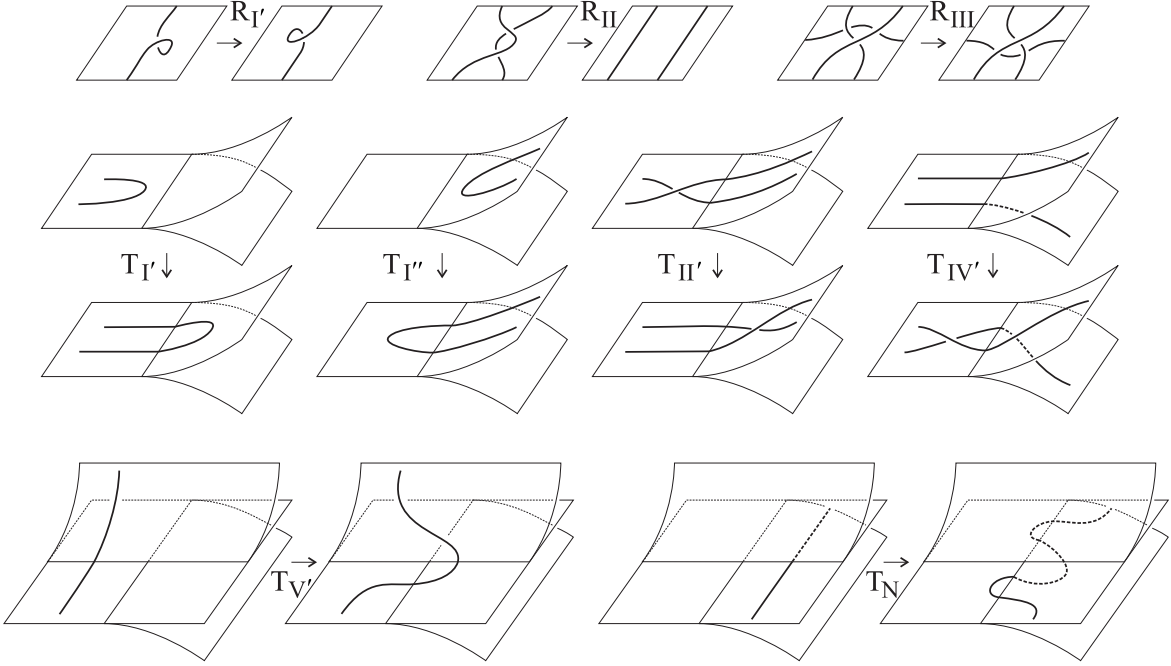,width=15.5cm}}}
\caption{\label{c1moves}$\cont^1$-Turaev moves.}
\end{figure}

\paragraph{Statement} 
Let us fix a branched standard
spine $P$ of $M$. The fact that such a spine always exists was proved as 
Theorem~3.4.9 in~\cite{lnm}. We call 
$\cont^1$ {\em link diagram on} $P$ an immersion of a disjoint union of
circles into $P$, with generic intersection with $S(P)$ appearing
as in Fig.~\ref{c1proj}, 
generic self-intersections (crossings), and the usual under-over marking at crossings
(as shown in the same figure; here `under' and `over' refer to the field positively
transversal to $P$).
The set of all $\cont^1$ link diagrams on $P$ will be denoted by $\dd(P)$.
An element of $\dd(P)$ obviously defines a link. Moreover $v(P)$ is transversal to 
this link, so it defines a framing, and we get an (obviously well-defined) map
$F_{(P,M)}:\dd(P)\to\Fram(M)$. 
Besides isotopy on $P$ through immersions having the same
configuration of crossings and intersections with $S(P)$, there are several combinatorial
moves which of course do not modify the isotopy class of the framed link defined by a 
diagram. We show a list of moves having this property in Fig.~\ref{c1moves}.
For a reason to be given below, which also explains the apparently weird notation, 
we call these moves $\cont^1$-Turaev moves.
As we did when we described the sliding-MP-moves in Fig.~\ref{slidingMP}, we are showing
in Fig.~\ref{c1moves} only the essential physical modifications, without specifying the 
screw-orientation of the spine and the orientation of its regions.

\begin{teo}\label{fram:statement}
The map $F_{(P,M)}:\dd(P)\to\Fram(M)$ is surjective, and the equivalence relation
it defines is generated by $\cont^1$-Turaev moves.
\end{teo}

Our argument, after the easy proof of surjectivity, goes along the following lines:
\begin{enumerate}
\item We state the analogue of Theorem~\ref{fram:statement} for non-branched spines,
due to Turaev~\cite{turaev:ombre} (we include a quick proof for the sake of completeness);
\item We modify Turaev's result to the case of a branched spine, but allowing 
non-$\cont^1$ diagrams;
\item We prove our theorem, showing how to canonically replace each non-$\cont^1$
diagram by a $\cont^1$ one along a sequence of modifications.
\end{enumerate}

\paragraph{Surjectivity} 
Since $M$ and $P$ are fixed, we write $F$ for short.
Given a framed link $L$ in $M$, we can prove that it is contained in the image of $F$
as follows:
\begin{itemize}
\item First, forget the framing and take a generic projection on $P$,
recalling that $M\setminus P\cong\partial M\times(0,1]$;
\item Next, eliminate non-$\cont^1$ intersections with $S(P)$ as 
shown in Fig.~\ref{goodproj} (left).
\begin{figure}
\figura{\vspace{1in}
\centerline{\special{wmf:SecondFigs/SecondWMF/goodproj.wmf y=1in}\hspace{5cm}\ }}
{\centerline{\psfig{file=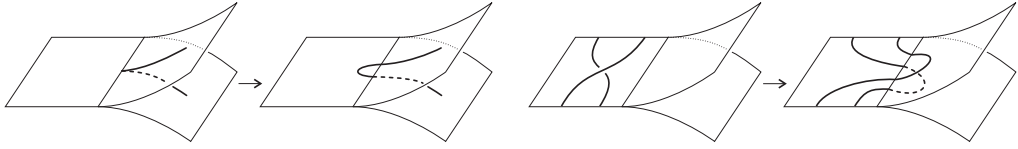,width=13.5cm}}}
\caption{\label{goodproj}How to remove forbidden intersections with 
$S(P)$ and crossings.}
\end{figure}
\item Finally, give the resulting projection the right framing by adding the
necessary numbers of curls. (Here we use the fact that two framings on a given
knot differ at most by a finite number of full rotations.)
\end{itemize}

It may be noted that surjectivity of $F$ is preserved by restriction 
to the set of diagrams without crossings. This follows quite easily from
the fact that all the regions of $P$ have non-empty boundary, as suggested in
Fig.~\ref{goodproj} (right). This property will not be used below.

\paragraph{Turaev moves on standard spines} 
The ideas and results of this paragraph
are due to Turaev~\cite{turaev:ombre}. We temporarily allow $P$
to be any standard spine of $M$, not a branched one. If each region of
$P$ is given an arbitrary transverse orientation, the definition of a
{\em link diagram} $D$ makes sense also on $P$, but $D$ may not
define a framing on the associated link, because the strip which runs along a component of $D$ 
on $P$ need not be a cylinder, it may be a M\"obius strip. So we 
attach to each component $D_i$ of
$D$ a full or half-integer $a_i/2$, depending on the topology of the strip,
and we define the framing by giving $a_i$ positive half-twists to the strip
(recall that $M$ is oriented). By {\em diagram} on $P$ we will
actually mean one such pair $(\{D_i\},\{a_i/2\})$.

\begin{figure}
\figura{\vspace{2in}
\centerline{\special{wmf:SecondFigs/SecondWMF/turmoves.wmf y=2in}\hspace{5cm}\ }}
{\centerline{\psfig{file=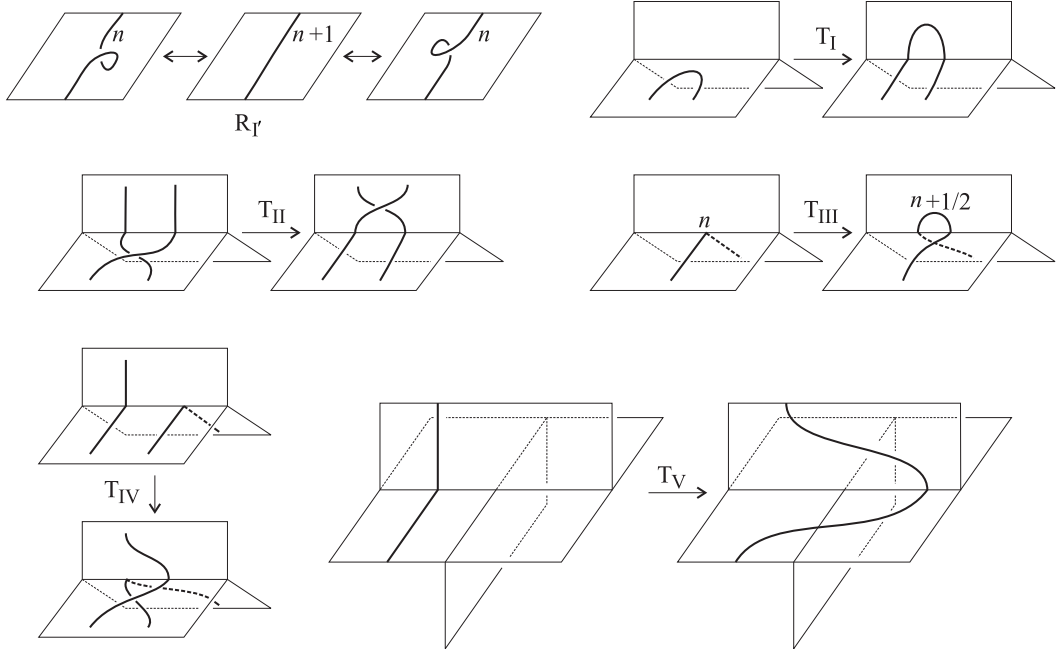,width=14cm}}}
\caption{\label{turmoves}Turaev moves.}
\end{figure}

We call Turaev moves those shown in Fig.~\ref{turmoves},
together with the ${\rm R}_{\rm I\!I}$ and ${\rm R}_{\rm I\!I\!I}$
already shown above. In Fig.~\ref{turmoves}, for ${\rm R}_{{\rm I}'}$ and
${\rm T}_{{\rm I\!I\!I}}$, the local orientation must be that of $\mr^3$.

\begin{teo}\label{tur:teo}
Every isotopy class of framed link in $M$ is defined by some diagram on $P$, and
two diagrams define the same class if and only if they are obtained from each
other by a sequence of Turaev moves.
\end{teo}

\begin{figure}
\figura{\vspace{.7in}
\centerline{\special{wmf:SecondFigs/SecondWMF/striproj.wmf y=.7in}\hspace{5cm}\ }}
{\centerline{\psfig{file=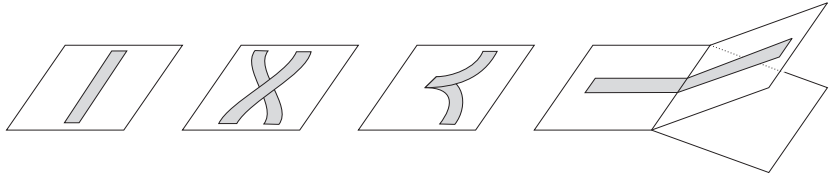,width=11cm}}}
\caption{\label{striproj}Projection of a strip.}
\end{figure}

\dim{tur:teo}
Recall that a framed link can be thought of as an embedded cylinder. Moreover $M$
projects onto $P$, and the projection of a cylinder generically appears as
in Fig.~\ref{striproj}.
Such a projection easily defines a diagram, the half-integers being sums of $\pm1/2$'s
corresponding to the bends of the projection. Moreover, 
the elementary catastrophes along an isotopy of a projection translate into the moves of Fig.~\ref{turmoves}, or simple combinations of them.\finedim{tur:teo}

\paragraph{Turaev moves on branched spines}
Going back to the case where $P$ is branched, we can still apply Theorem~\ref{tur:teo},
but now the list of moves becomes slightly longer, if we want to take the branching into
account.

\begin{prop}\label{short:list}
If $P$ is a branched spine then any Turaev move for a diagram
on $P$ can be expressed as a combination (including inverses)
of the moves ${\rm R}_{{\rm I}'}$,
${\rm R}_{\rm I\!I}$, ${\rm R}_{\rm I\!I\!I}$,
${\rm T}_{{\rm I}'}$, ${\rm T}_{{\rm I}''}$, ${\rm T}_{{\rm I\!I}'}$,
${\rm T}_{{\rm I\!V}'}$, ${\rm T}_{{\rm V}'}$
shown above, together with the moves 
${\rm T}_{{\rm I}'''}$, ${\rm T}_{{\rm I\!I}''}$, ${\rm T}_{{\rm I\!I\!I}'}$, 
${\rm T}_{{\rm I\!V}''}$, ${\rm T}_{{\rm V}''}$
shown in Fig.~\ref{mixmoves}
\begin{figure}
\figura{\vspace{2in}
\centerline{\special{wmf:SecondFigs/SecondWMF/mixmoves.wmf y=2in}\hspace{5cm}\ }}
{\centerline{\psfig{file=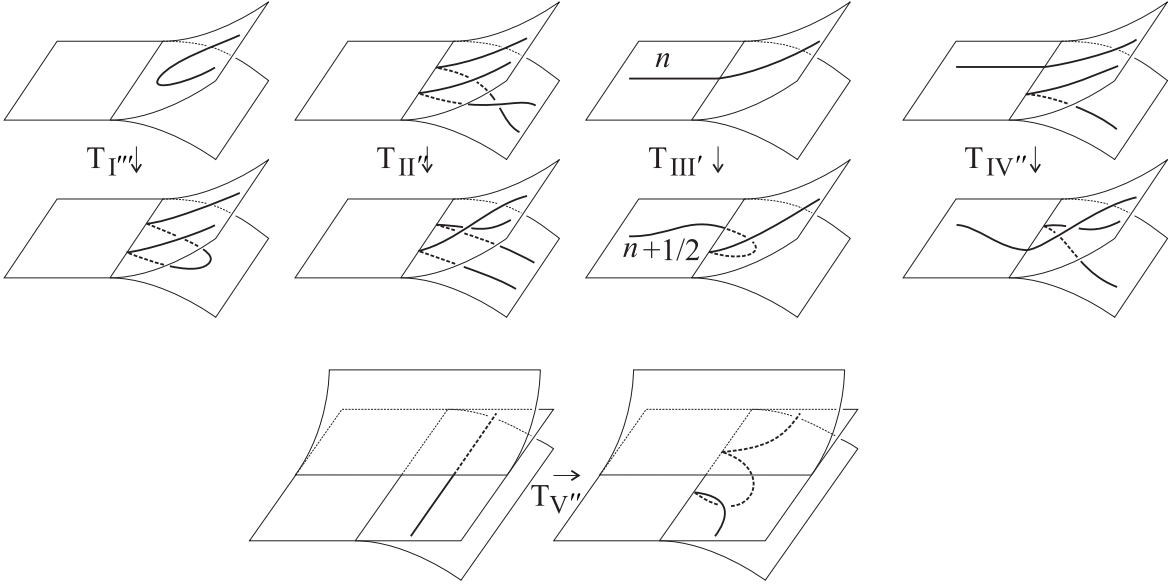,width=15.5cm}}}
\caption{\label{mixmoves}Non-$\cont^1$ Turaev moves on branched spines.}
\end{figure}
\end{prop}

\dim{short:list}
The branching can be interpreted as a loss of symmetry of a spine, so
each of Turaev's moves, when viewed as a move on a branched spine,
generates many different ones according to the position of the diagram
with respect to the branching. The result is a list much longer than that given
in the statement, but one can show that all the moves omitted from the statement
\begin{figure}
\figura{\vspace{.7in}
\centerline{\special{wmf:SecondFigs/SecondWMF/omitex1.wmf y=.7in}\hspace{5cm}\ }}
{\centerline{\psfig{file=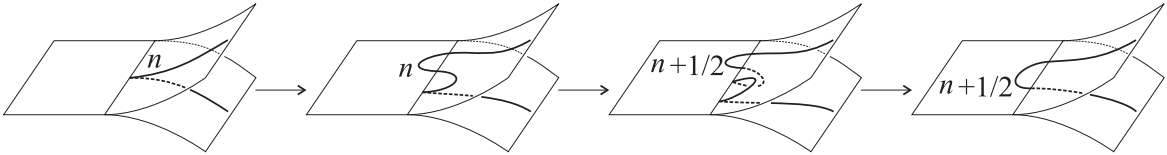,width=15.5cm}}}
\caption{\label{omitex1}A sequence of moves ${\rm T}_{{\rm I}''}$, 
${\rm T}_{{\rm I\!I\!I}'}$ and ${\rm T}_{{\rm I}'''}$ generating a branched version of
${\rm T}_{\rm I\!I\!I}$ omitted from the statement.}
\end{figure}
are generated by the moves included. Two examples are provided in Figg.~\ref{omitex1}
\begin{figure}
\figura{\vspace{.7in}
\centerline{\special{wmf:SecondFigs/SecondWMF/omitex2.wmf y=.7in}\hspace{5cm}\ }}
{\centerline{\psfig{file=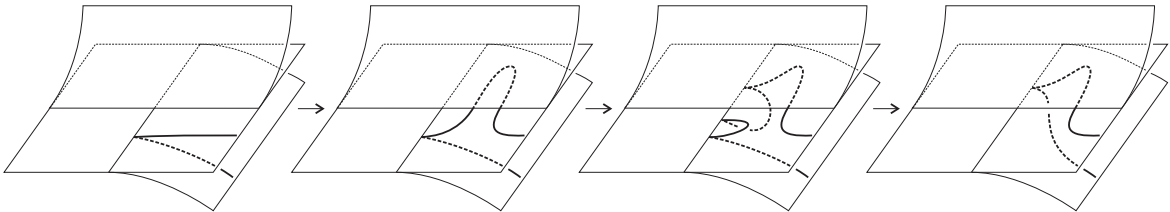,width=15.5cm}}}
\caption{\label{omitex2}Another generation of moves, involving ${\rm T}_{{\rm I}'}$, 
${\rm T}_{{\rm V}''}$ and ${\rm T}_{{\rm I}'''}$.}
\end{figure}
and~\ref{omitex2}.
\finedim{short:list}

\paragraph{${\bf C}^{\bf 1}$ moves} 
We can now conclude the proof of Theorem~\ref{fram:statement}
(surjectivity having been shown above). We are left to show that if $D$ and $D'$ are
$\cont^1$-diagrams on $P$ which define the same framed link, then they
are related by a sequence of $\cont^1$ Turaev moves. By Theorem~\ref{tur:teo}
and Proposition~\ref{short:list}, there exists a sequence
$D=D_0\to D_1\to\cdots\to D_{n-1}\to D_n=D'$ where each move $D_{i-i}\to D_i$ is one
of those listed in Proposition~\ref{short:list}. In particular the $D_i$'s with $0<i<n$
can be non-$\cont^1$ and can have a non-zero half-integer attached to them.
We will now show how to construct
a modified sequence $D=\tilde D_0\to\tilde D_1\to\cdots\to\tilde D_{n-1}\to \tilde D_n$ 
with the following properties:
\begin{enumerate}
\item each $\tilde D_i$ is a $\cont^1$ diagram with number 0 attached;
\item each $\tilde D_i$ is obtained from $\tilde D_{i-1}$ by a sequence of
$\cont^1$ Turaev moves;
\item\label{pre-key:step} each $\tilde D_i$ differs from 
$D_i$ for the presence of some extra curls;
in particular each component $\tilde D^{(j)}_i$ of $\tilde D_i$ has a natural companion
$D^{(j)}_i$ in $D_i$, with the property that, as unframed knots, the knots associated to 
$\tilde D^{(j)}_i$ and $D^{(j)}_i$ are both contained in a solid torus $T^{(j)}_i$
and parallel to the core of the torus;
\item\label{key:step} the framed knots associated to $\tilde D^{(j)}_i$ and $D^{(j)}_i$ 
{\em are framed-isotopic within $T^{(j)}$}.
\end{enumerate}
Requirement~\ref{key:step} is the crucial technical point of our proof. To verify that 
the requirement is stronger than just framed isotopy,
note that in $D^2\times S^1$ the framings on the core $\{0\}\times S^1$ are
parameterized by the integers, but, when $D^2\times S^1$
is mirrored in its boundary to get $S^2\times S^1$,
only two inequivalent framings remain (corresponding to even and odd integers).

We assume for a moment the sequence $\tilde D_i$ to exist, and we show
how to conclude. The transformation from $D$ to $\tilde D_n$ is made with $\cont^1$
Turaev moves, so we only need to compare $\tilde D_n$ and $D'=D_n$, which by assumption
differ for some curls. Using $\cont^1$ Turaev moves 
we can easily make all these curls slide until they are consecutive on the
diagram. We recall
now that there are four local pictures for a curl, depending on its
local contributions $\pm$ to the framing and to the winding number~\cite{trace}.
Assumption~\ref{key:step} now implies that the algebraic sum of local contributions
to the framing vanishes. Therefore we can cancel the curls in pairs, either by
moves ${\rm R}_{{\rm I}'}$ (when the local contributions to the winding number
are the same), or by a combination of moves ${\rm R}_{\rm I\!I}$ and
${\rm R}_{\rm I\!I\!I}$ (when the contributions cancel). This shows the conclusion.

We are left to define the sequence $\tilde D_i$. The idea is simply not to
perform the moves which change the half-integer or introduce 
cusps, and show that the sequence of moves can be followed anyway.
While doing this we need to keep track of the portions where the new
diagram $\tilde D_i$ differs from $D_i$, which we do by marking a neighbourhood 
of the portion as a shadowed box.

\begin{figure}
\figura{\vspace{.7in}
\centerline{\special{wmf:SecondFigs/SecondWMF/replam1.wmf y=.7in}\hspace{5cm}\ }}
{\centerline{\psfig{file=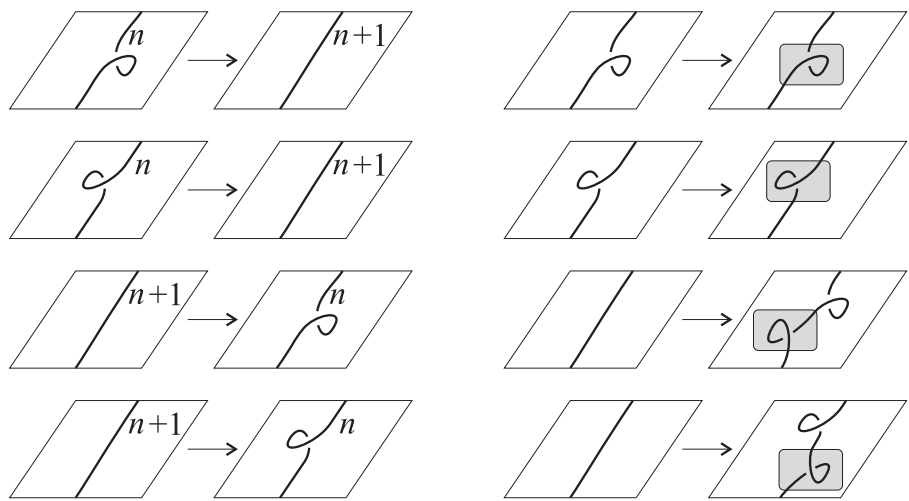,width=12cm}}}
\caption{\label{replam1}The moves ${\rm R}_{{\rm I}}$ (left), and the $\cont^1$ moves
replacing them (right).}
\end{figure}

\begin{figure}
\figura{\vspace{.7in}
\centerline{\special{wmf:SecondFigs/SecondWMF/replam2.wmf y=.7in}\hspace{5cm}}}
\centerline{\psfig{file=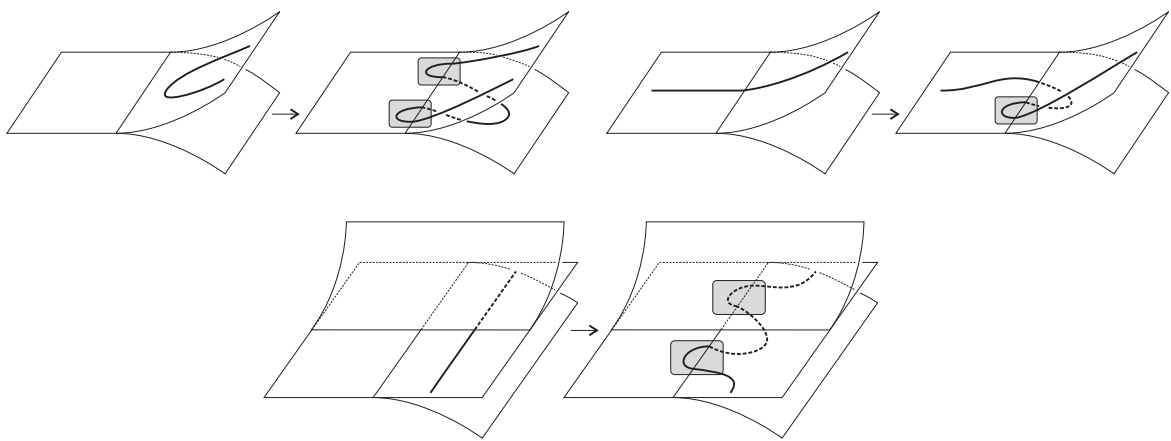,width=15.5cm}}
\caption{\label{replam2}Replacement of the other non-$\cont^1$ moves ${\rm T}_{{\rm I}'''}$, ${\rm T}_{{\rm I\!I\!I}'}$ and ${\rm T}_{{\rm V}''}$ by $\cont^1$ moves.}
\end{figure}

The moves which change the colour or introduce cusps are ${\rm R}_{{\rm I}'}$
(in both directions), ${\rm T}_{{\rm I}'''}$, ${\rm T}_{{\rm I\!I\!I}'}$,
and ${\rm T}_{{\rm V}''}$, and we show in Figg.~\ref{replam1} and~\ref{replam2} 
what we replace them with. For move ${\rm R}_{{\rm I}'}$ it has been necessary to be more
specific because, in the original definition of the move, two different ones
were actually defined at the same time.

To show that $\{\tilde D_i\}$ can indeed be constructed we must now show that
after performing the construction up to some level $k$ we can still 
still follow the rest of the sequence and go on with the replacements of
Figg.~\ref{replam1} and~\ref{replam2}. By construction $\tilde D_k$ differs from
$D_k$ only within some shadowed boxes. We denote by $\mu_k$ the move
$D_k\to D_{k+1}$ and explain how to lift it to a move on $\tilde D_k$. First,
note that a shadowed box lying
within a region of $P$ and containing a curl
does not interfere with $\mu_k$ whatever its type (but it may be necessary
to add some ${\rm R}_{\rm I\!I}$'s and ${\rm R}_{\rm I\!I\!I}$'s to replace isotopy
supported within the region). 

We fix now our attention on a shadowed box $B$ which lies
on $S(P)$ and contains a smoothed cusp, and examine the various instances for
$\mu_k$ with respect to $B$. If $\mu_k$ is of type 
${\rm R}_{\rm I\!I}^{\pm1}$, ${\rm R}_{\rm I\!I\!I}^{\pm1}$,
${\rm T}_{{\rm I}'}^{\pm1}$, ${\rm T}_{{\rm I}''}^{\pm1}$
${\rm T}_{{\rm I\!I}'}^{\pm1}$, or ${\rm T}_{{\rm I\!V}'}^{\pm1}$, 
then it obviously does not interfere with $B$, so we can just
perform $\mu_k$ on $\tilde D_k$. 
If $\mu_k$ is one of the moves ${\rm R}_{\rm I}^{\pm1}$, ${\rm T}_{{\rm I}'''}$,
${\rm T}_{{\rm I\!I\!I}'}$ or ${\rm T}_{{\rm V}''}$ then again it does not interfere
with $B$, and we can perform the appropriate replacement from Figg.~\ref{replam1}
or~\ref{replam2}, getting a move from $\tilde D_k$ to $\tilde D_{k+1}$.
If $\mu_k$ is a move of type
${\rm T}_{{\rm I\!I}''}^{\pm1}$ or ${\rm T}_{{\rm I\!V}''}^{\pm1}$
then it may interfere with $B$. However, since it does not create or destroy
cusps, $\mu_k$ can be translated on $\tilde D_k$ as a combination of moves 
which do not involve cusps, {\em i.e.} allowed from the statement.
We are only left to deal with the case where $\mu_k$ is one of the moves
${\rm T}_{{\rm I}'''}^{-1}$,
${\rm T}_{{\rm I\!I\!I}'}^{-1}$ or ${\rm T}_{{\rm V}''}^{-1}$ which destroy cusps.
By construction the cusp(s) to be destroyed still appear in $\tilde D_k$ as
smoothed cusps within shadowed boxes, so we can destroy them also from
$\tilde D_k$ by means of allowed moves.

The crucial properties~\ref{pre-key:step} and~\ref{key:step} hold by construction,
and the conclusion eventually follows.\finedim{fram:statement}

\section{A calculus for pseudo-Legendrian links\\
in combed manifolds}\label{leg:calc}
We will deal in this section with the set $\Pleg$ of 
equivalence classes of triples $(M,v,L)$ already described in the introduction.
Its combinatorial counterpart will be given by the set
$$\call=\{(P,D):\ P\in\bb,\ D\in\dd(P)\}$$
where $\bb$ is as in Section~\ref{conc:calc} and $\dd(P)$ is as  in
Section~\ref{fram:calc}. The reconstruction map $(P,D)\mapsto L(P,D)$
is here defined by noting that $D$ defines a link transversal to $v(P)$
in $M(P)$, and hence also a link transversal to $\hatv(P)$ in $\hatM(P)$.
According to what we stated after Proposition~\ref{rigid:spines}
we will actually drop from $\bb$ the sliding-MP-rigid spines, and 
ignore the snake move.

If for a fixed $P$ we consider the effect on $L(P,D)$ of the 
$\cont^1$-Turaev moves on $D$, we see that the class 
of $L(P,D)$ is in general modified by the 
first Reidemeister move ${\rm  R}_{{\rm I}'}$ of Fig.~\ref{c1moves}
(see also Section~\ref{specu:section}), but not
by the other moves, which we will therefore call {\em pseudo-Legendrian Turaev moves}.
Other moves which obviously do not change $L(P,D)$ up to equivalence are
the sliding-MP-moves on $P$ which do not
involve $D$ (these moves permit to follow $D$ along the modification of $P$, so
they are well-defined for pairs). It is not hard to see that 
before performing a sliding-MP-move on $P$ it is always possible to modify $D$
by pseudo-Legendrian Turaev moves to a diagram which is not involved in the sliding-MP-move,
and the diagram after the sliding-MP-move is well-defined up to 
pseudo-Legendrian Turaev moves.
For this reason we will freely speak of sliding-MP-moves also for pairs.
The following will be established below.

\begin{teo}\label{leg:calc:teo}
The map $L:\call\to\Pleg$ is surjective, and the equivalence relation
defined by $L$ is generated by pseudo-Legendrian Turaev moves and sliding-MP-moves.
\end{teo}

\paragraph{Fixed-spine statement} 
Recall from the definition that the 
map $L$ of the statement of Theorem~\ref{leg:calc:teo} involves the passage from
$P$ to $(M(P),v(P))$ and then to $(\hatM(P),\hatv(P))$. As already pointed out,
this is necessary if one wants to be able to deal with non-traversing fields.
However, if one happens to have a concave traversing field, one can directly encode
this field by a spine, without first removing a ball, and one can investigate
how isotopy of links transversal to the field reflects on link diagrams on the
spine. The following is shown below:

\begin{prop}\label{fixed:P:prop}
Let $P$ be a branched spine. Fix a representative of $(M(P),v(P))$ and
an embedding of $P$ in $M(P)$ transversal to $v(P)$. Then every link transversal
to $v(P)$ is represented by a $\cont^1$ diagram on $P$. Moreover two 
$\cont^1$ diagrams define the same link up to isotopy through links transversal
to $v(P)$ if and only if they are related by pseudo-Legendrian Turaev moves.
\end{prop}

\paragraph{From fixed to variable spine} 
We show in this
paragraph how to deduce Theorem~\ref{leg:calc:teo} from
Proposition~\ref{fixed:P:prop}. First of all, to prove surjectivity, we consider
a triple $(M,v,L)$ representing an element of $\Pleg$.
Using a normal section as in Proposition~\ref{section:then:spine}, we can obtain
a spine $P\in\bb$ which encodes the equivalence class of $(M,v)$ in the sense of Theorem~\ref{comb:calc:teo}. Moreover $P$ comes with an embedding
in $M$ transversal to $v$. Now, a neighbourhood of $P$ can be identified to
$M(P)$ and its complement is isomorphic to $\bthtr$. Using the flow generated by $v$ in this
ball we can now isotope $L$ through links transversal to $v$ to a link
which lies in $M(P)$, and the first assertion of Proposition~\ref{fixed:P:prop}
implies that $L$ is represented by a diagram $D$ on $P$. Summing up, we see that
$(M,v,L)$ is represented by $(P,D)$, and surjectivity of $L$ is proved.

To conclude we must now show that two pairs $(P_0,D_0)$ and $(P_1,D_1)$
are equivalent via pseudo-Legendrian Turaev moves when $L(P_0,D_0)=L(P_1,D_1)$.
Spelling out the relation of pseudo-Legendrian isotopy, which defines $\Pleg$,
we assume that $P_0$ and $P_1$ embed in the same manifold $M$ and that there exist
a field homotopy $(v_t)_{t\in[0,1]}$ and a link isotopy $(L_t)_{t\in[0,1]}$ on $M$
such that:
\begin{enumerate}
\item for $i=0,1$, the link $L_i$ is the one defined by $D_i$, and the field
$v_i$ is positively transversal to
$P_i$ and restricts to $\bthtr$ on the complement of $P_i$;
\item $L_t$ is transversal to
$v_t$ for all $t$.
\end{enumerate}
Given $t\in[0,1]$, we note that for $|t-s|\ll 1$ the link 
$L_s$ is transversal to $v_t$, and that a branched spine for $v_t$ 
(in the sense repeatedly used above) is a branched spine also for $v_s$. So we can
subdivide $[0,1]$ into subintervals $[t_{i-1},t_i]$ so that:
\begin{enumerate}
\item $L_{t_{i-1}}$ is isotopic to $L_{t_i}$ through links transversal to $v_{t_i}$;
\item $v_{t_i}$ has a spine $P_{t_i}$ which is also a spine for $v_{t_{i-1}}$.
\end{enumerate}
Now let $D_{t_i}$ be a diagram for $L_{t_i}$ on $P_{t_i}$. Since both $P_{t_{i-1}}$ and
$P_{t_i}$ are spines for $v_{t_i}$, we can transform $P_{t_{i-1}}$ into 
$P_{t_i}$ via sliding-MP-moves. Using pseudo-Legendrian Turaev moves we can now
follow $D_{t_{i-1}}$ along this sequence of moves, getting a diagram $D'_{t_i}$ on $P_{t_i}$.
Since the sequence of sliding-MP-moves can be realized in $M$ so that each spine
of the sequence is a branched spine for $v_{t_i}$, we deduce that 
the link defined by $D'_{t_i}$ is isotopic to $L_{t_{i-1}}$, and hence to $L_{t_i}$,
through links transversal to $v_{t_i}$. Proposition~\ref{fixed:P:prop}
now implies that $D'_{t_i}$ and $D_{t_i}$ are related by pseudo-Legendrian Turaev moves on
$P_{t_i}$. This shows that $(P_{t_i},D_{t_i})$ 
is obtained from $(P_{t_{i-1}},D_{t_{i-1}})$ via the moves of the statement,
and the conclusion follows by iteration.\finedim{leg:calc:teo}

\paragraph{Fixed-spine proof} 
We will establish now Proposition~\ref{fixed:P:prop},
writing just $M$ and $v$ for $M(P)$ and $v(P)$.
For the sake of simplicity we will assume that $v$ is tangent to $\partial M$
along only one curve (denoted by $\gamma$), but our arguments 
extends almost {\em verbatim} to the general case of more than one curve. 

\begin{figure}
\figura{\vspace{.7in}
\centerline{\special{wmf:SecondFigs/SecondWMF/annulus.wmf y=.7in}\hspace{5cm}\ }}
{\centerline{\psfig{file=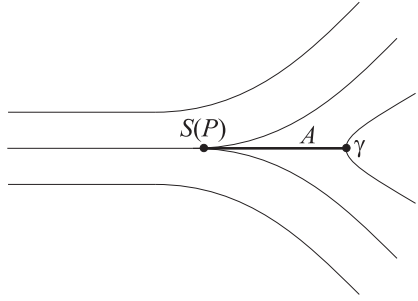,width=5.5cm}}}
\caption{\label{annulus}The annulus $A$.}
\end{figure}

We fix in $M$ an annulus $A$ which connects $\gamma$ to $S(P)$ as shown
in a cross-section in Fig.~\ref{annulus}.
Note that $A$ is almost but not quite embedded: it has double point at the vertices of $P$.
Since we will only need to consider $A$ locally and away from vertices of $P$, this
fact will not disturb us. We choose coordinates $(\rho,\theta)\in[0,1]\times[0,2\pi]$
on $A$, where $\rho=0$ corresponds to $S(P)$ and $\rho=1$ to $\gamma$. Near $A$ we can also define a coordinate $z\in[-\varepsilon,\varepsilon]$ by integrating $v$.

\begin{figure}
\figura{\vspace{.7in}
\centerline{\special{wmf:SecondFigs/SecondWMF/al_types.wmf y=.7in}\hspace{5cm}\ }}
{\centerline{\psfig{file=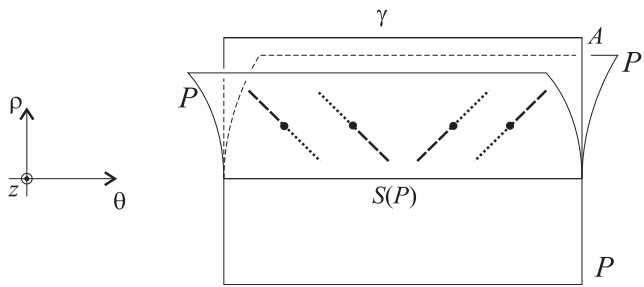,width=8.5cm}}}
\caption{\label{al_types}Types of points in $L\cap A$.}
\end{figure}

\begin{figure}
\figura{\vspace{.7in}
\centerline{\special{wmf:SecondFigs/SecondWMF/par_proj.wmf y=.7in}\hspace{5cm}\ }}
{\centerline{\psfig{file=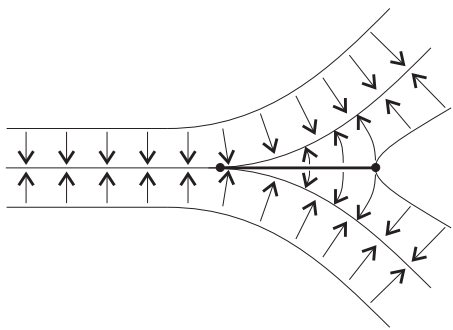,width=6cm}}}
\caption{\label{par_proj}Projection onto $P$.}
\end{figure}

\begin{figure}
\figura{\vspace{1.5cm}
\hspace{-4cm}
\centerline{\special{wmf:SecondFigs/SecondWMF/specproj.wmf y=1.5cm}
\hspace{10cm}
\special{wmf:SecondFigs/SecondWMF/geneproj.wmf y=1.5cm}
\hspace{.5cm}\ }}
{\centerline{\psfig{file=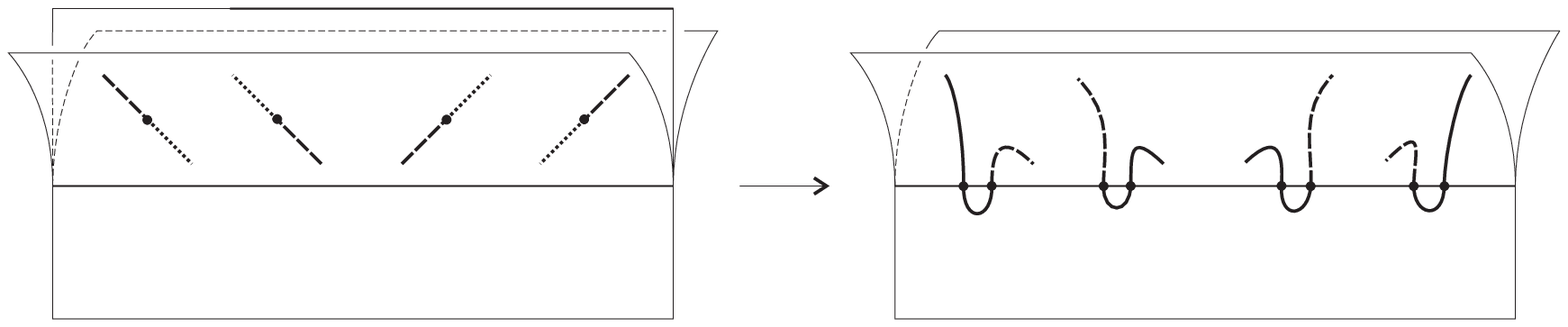,height=2.5cm}}
\vspace{.5cm}
\centerline{\psfig{file=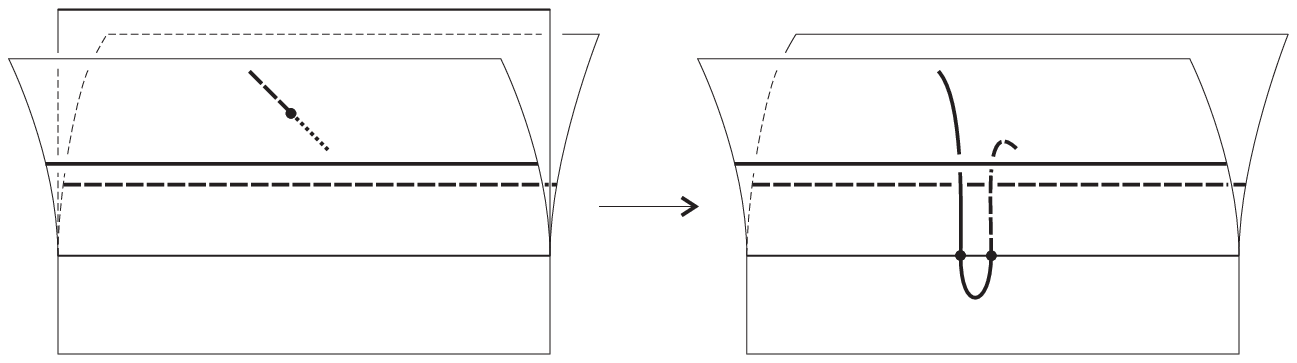,height=2.5cm}}}
\caption{\label{specproj}Completion of the diagram at points of $L\cap A$.}
\end{figure}

\begin{figure}
\figura{\vspace{.7in}
\centerline{\special{wmf:SecondFigs/SecondWMF/accident.wmf y=.7in}\hspace{5cm}\ }}
{\centerline{\psfig{file=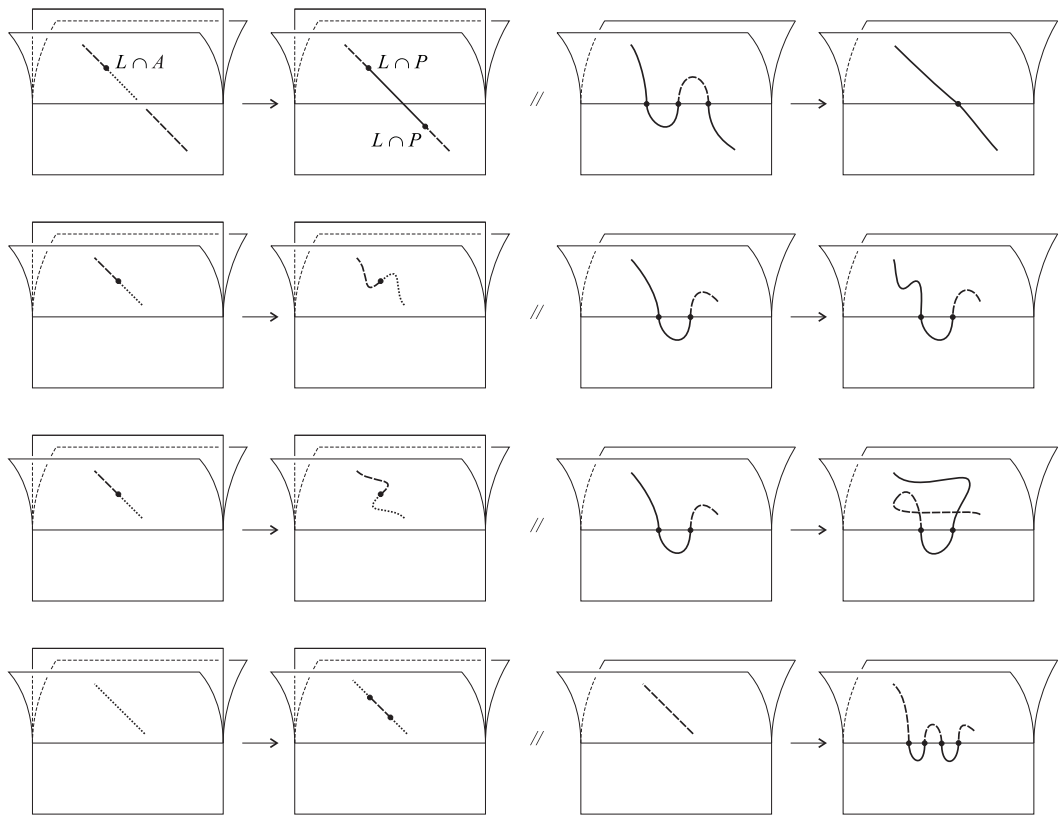,width=14cm}}}
\caption{\label{accident}Elementary accidents involving $A$ along an isotopy of $L$.}
\end{figure}

Now let $L$ be transversal to $v$, and assume by general position that $L$ intersects $A$
only at points with $0<\rho<1$, that no two such intersections have the same coordinate
$\theta$, and that at all the intersections the tangent direction to $L$ has
non-zero components in all three coordinates $\rho,\theta,z$.
Depending on the sign of these components, we can divide the points of $L\cap A$
into four types, shown in Fig.~\ref{al_types}
($L$ is dashed when it lies over $A$ and dotted when it lies under $A$).
We consider now the projection $\pi$ of $M\setminus A$ onto $P$ along the orbits of $v$,
as shown in Fig.~\ref{par_proj}.
Of course $L\setminus A$ locally projects to a $\cont^1$-strand on $P$, and
by general position we can assume that $\pi(L\setminus A)$ locally appears as a 
$\cont^1$-diagram. We are only left to extend the diagram at the points of $L\cap A$, 
which we do locally in Fig.~\ref{specproj}.
The top part of
this figure actually refers to a simplified situation, because other strands of $L$
already projected on $P$ may locally interfere. We show at the
bottom of the same figure
in one example how to deal with this fact.
The resulting diagram of course represents $L$, and we have proved the first
assertion in Proposition~\ref{fixed:P:prop}. To prove the second assertion we
must now examine an isotopy of $L$ through links transversal to $v$, and
hence examine first-order violations of genericity of $L$ with respect to $A$
and $\pi$. All the elementary accidents which do not involve $A$ of course
correspond to pseudo-Legendrian Turaev moves. We are left to deal with the following accidents:
\begin{enumerate}
\item $L$ intersects $A$ at a point of $S(P)\subset\partial A$;
\item at a point of $L\cap A$, the tangent direction to $L$
has vanishing $\rho$-coordinate;
\item similarly, with the $\theta$-coordinate;
\item similarly, with the $z$-coordinate.
\end{enumerate}
In Fig.~\ref{accident} we show the situation
just before and just after each of these accidents, and
we analyze the corresponding transformations of the diagrams constructed 
as in Fig.~\ref{specproj}. In all cases one easily sees that indeed the 
transformation is generated by pseudo-Legendrian Turaev moves: the number of moves
needed is respectively one, zero (isotopy within regions), three and two. By simplicity in
Fig.~\ref{accident} we have ignored the possible interference of other strands of $L$ 
projected on $P$, but the conclusion is valid anyway (some Reidemeister moves must be added
in the general case).\finedim{fixed:P:prop}

\section{Applications and speculations}\label{specu:section}
In this section we will discuss some consequences of the calculi
described above, and mention some natural questions and problems which
we put forward for further investigation. The section is split into two subsections.

\subsection{Winding number, torsion, and \foi s}\label{wind:section}
In this section we employ our realizations of $\Fram$ and $\Pleg$
in connection with winding number, Maslov index, torsion, and \foi s of
pseudo-Legendrian knots.

\paragraph{Relative winding number}
We spell out in this paragraph the analogue of Trace's result~\cite{trace}
on knot diagrams in $\mr^3$. We confine ourselves to knots for the
sake of simplicity, but essentially the same holds for links.

\begin{prop}\label{trace:like}
Let $(v_0,K_0)$ and $(v_1,K_1)$ be pseudo-Legendrian pairs in a manifold $M$,
where $K_0$ and $K_1$ are oriented knots. Assume that $v_0$ and $v_1$ are homotopic
relatively to $\partial M$, and that $K_0$ and $K_1$ are isotopic as oriented
framed knots. Then, up to pseudo-Legendrian isotopy on $(v_1,K_1)$, we can
assume that $v_1=v_0$ and that $K_1$ differs from $K_0$ only within a region
of $(M,v_0)$ isomorphic to $(\mr^3,\partial/\partial z)$, where $K_0$ is a straight
horizontal line and $K_1$ has either some positive or some negative double curls
(shown in Fig.~\ref{twocurls}).
\begin{figure}
\figura{\vspace{.2in}
\centerline{\special{wmf:SecondFigs/SecondWMF/twocurls.wmf y=.2in}\hspace{5cm}\ }}
{\centerline{\psfig{file=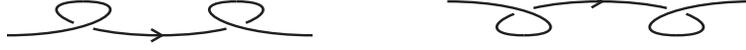,width=10cm}}}
\caption{\label{twocurls}A positive and a negative double curl (here the field
points upwards and the orientation is consistent with the orientation of 3-space).}
\end{figure}
\end{prop}

\begin{figure}
\figura{\vspace{.3in}
\centerline{\special{wmf:SecondFigs/SecondWMF/avoidRI.wmf y=.3in}\hspace{5cm}\ }}
{\centerline{\psfig{file=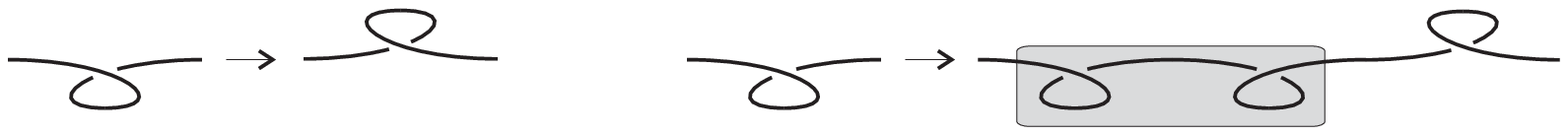,width=14cm}}}
\caption{\label{avoidRI}Replacing a non-pseudo-Legendrian move.}
\end{figure}

\dim{trace:like}
Choose branched spines $P_0$ and $P_1$ of $v_0$ and $v_1$ according to 
Proposition~\ref{section:then:spine}, and use Proposition~\ref{fixed:P:prop} to
represent $K_0$ and $K_1$
on $P_0$ and $P_1$ respectively by $\cont^1$-diagrams $D_0$ and $D_1$.
Since $v_0$ and $v_1$ are homotopic, a sequence of sliding-MP moves
connects $P_1$ to $P_0$. Following $D_1$ along this sequence of moves
we get a pseudo-Legendrian isotopy, so we can assume that $v_1=v_0$ and $P_1=P_0$.
Now $D_0$ and $D_1$ define on $P_0$ framed-isotopic oriented knots, so
by Theorem~\ref{fram:statement} they are related by $\cont^1$-Turaev moves.
If along the sequence of moves there is no ${\rm R}_{{\rm I}'}$, we deduce
pseudo-Legendrian equivalence of $K_0$ and $K_1$. If however there is some
${\rm R}_{{\rm I}'}$, we replace it by a pseudo-Legendrian move as shown in
Fig.~\ref{avoidRI}.
One easily sees that this replacement can be done consistently along the sequence
of moves. The result is a pseudo-Legendrian isotopy between $D_1$ and a diagram
which differs from $D_0$ only for some double curls. These double curls can 
of course be slid to be consecutive. Now, it is precisely the content of~\cite{trace}
that up to moves ${\rm R}_{\rm I\!I}$ and ${\rm R}_{\rm I\!I\!I}$
there are only the types of double curls shown in Fig.~\ref{twocurls},
and that a positive and a negative double curl cancel out.\finedim{trace:like}

Under the assumptions of the previous proposition one may be tempted to define
a relative winding number $w(K_1,K_0)$ as the (algebraic) number of double curls by 
which the diagram of $K_1$ differs from the diagram of $K_0$. This number is however not 
well-defined in general, as one easily sees in $S^2\times S^1$ with vector field
parallel to the $S^1$-factor, because in this case
a double curl on a diagram contained in $\mr^2=S^2\setminus\{\infty\}$ can
always be removed by isotoping the diagram through $\infty$. This seems to suggest
that not $w(K_1,K_0)$, but maybe 
$w(K_1,K_0)\cdot[\mu_{K_0}]\in H_1(E(K_0);\mz)$ 
is well-defined, where $E(K_0)$ is the exterior of $K_0$ and $\mu_{K_0}$ is the meridian.
We will show this fact under the additional assumption that $K_0$ is `good'
(see~\cite{first:paper} and below for explanations).

\paragraph{Torsion invariants and good knots}
In~\cite{first:paper} we have defined the 
Reidemeister-Turaev torsion of an Euler structure
with simple boundary, and we have applied this notion to define the torsion of
pseudo-Legendrian knots. As an absolute invariant torsion
contains a sort of lift of the classical Alexander invariant.
We will discuss in this paragraph the information carried by torsion
as a relative invariant of two pseudo-Legendrian
framed-isotopic knots $K_0,K_1$ in the same concave combed manifold $(M,v)$.
We recall from~\cite{first:paper} that this information is most easily expressed
when $K_0$ has the property of being {\it good}. Goodness depends only on the
isotopy class of $K_0$ as a framed knot, and it means that a certain quotient of
the mapping class group of $E(K_0)$ acts trivially on the space of Euler structures on
$E(K_0)$. We omit the precise definition here, but we recall that many knots
indeed are good (for instance, all are good if $M$ is a homology sphere, and
most hyperbolic knots are good).

When $K_0$ is good, the information carried by torsion as a relative invariant
depends only on $\alpha(v\ristr{E(K_0)},f(v\ristr{E(K_1)})\in H_1(E(K_0);\mz)$, 
where $f\in{\rm Diff}_0(M)$ maps $K_1$ to $K_0$ as framed knots, and
$\alpha$ is the first obstruction for two vector fields to be homotopic relative to the
boundary. So the next result means that for good knots the relative 
winding number gives a well-defined invariant, and all
the information torsion can capture is contained in the relative winding number.
The statement involves all the assumptions and notations of the present and previous
paragraph.

\begin{prop}\label{comp:class:prop}
$\alpha(v\ristr{E(K_0)},f(v\ristr{E(K_1)})=
w(K_1,K_0)\cdot[\mu_{K_0}]\in H_1(E(K_0);\mz).$
\end{prop}

\dim{comp:class:prop}
We note first that $w(K_1,K_0)$ and $[\mu_{K_0}]$ depend on the choice of an orientation
on $K_0$, but their product does not, so the statement makes sense. For the proof,
note that by goodness we can just assume that $K_0$ and $K_1$ differ as in the statement
of Proposition~\ref{trace:like}, and that $f$ is supported on a neighbourhood of the
region where $K_0$ and $K_1$ differ. The conclusion then follows directly from
Proposition~2.17 of~\cite{first:paper}.\finedim{comp:class:prop}

Even if we have not discussed \foi s yet, we note here that Proposition~\ref{comp:class:prop}
implies that torsion is a weaker 
invariant than the finite-order ones for Legendrian knots in a given
homotopy class of Legendrian immersions.
To our knowledge the only known examples of framed-isotopic knots
distinguished by such invariants are those due to Tchernov~\cite{vlad}, and we believe that 
they are all good (at least, they certainly are good when the ambient manifold is
$S^2\times S^1$). Now 
one sees that in all of Tchernov examples $w(K_1,K_0)\cdot[\mu_{K_0}]=0$, so torsion 
definitely cannot distinguish. On the other hand, the definition of
torsion does not require fixing a homotopy class of Legendrian immersions, so
torsion and \foi s are in some sense complementary.

We will state in the rest of this paragraph 
some interesting consequences of Proposition~\ref{comp:class:prop},
always assuming the knots involved to be good. For simplicity,
as in Proposition~\ref{comp:class:prop},
we stick to knots transverse to a given field $v$ on a given $M$,
but we remind that the relation of pseudo-Legendrian isotopy also involves
a homotopy of $v$.

\begin{cor}\label{zero:wind}
Under the same assumptions as in Proposition~\ref{comp:class:prop},
suppose furthermore that $[\mu_{K_0}]$ has infinite order in $H_1(E(K_0);\mz)$,
so $w(K_1,K_0)\in\mz$ is well-defined. Then the following facts are pairwise equivalent:
\begin{enumerate}
\item\label{zero:w:point} $w(K_1,K_0)=0$;
\item\label{zero:torsion:point} $K_0$ and $K_1$ have trivial relative torsion invariants;
\item\label{iso:point} $K_0$ and $K_1$ are pseudo-Legendrian isotopic.
\end{enumerate}
\end{cor}

\dim{zero:wind}
Equivalence of (\ref{zero:w:point}) and 
(\ref{iso:point}) follows from the definition of $w(K_1,K_0)$. 
Implication (\ref{zero:w:point})$\Rightarrow$(\ref{zero:torsion:point})
follows from Proposition~\ref{comp:class:prop}, and the
opposite implication follows by taking the torsion associated
to a representation $\varphi$ of $H_1(E(K_0);\mz)$ such that
$\varphi([\mu_{K_0}])$ has infinite order (see~\cite{first:paper} 
for details).\finedim{zero:wind}

If $M$ is a homology sphere and $K$ is a
pseudo-Legendrian knot in $(M,v)$ we have shown in~\cite{first:paper} that
the rotation number ${\rm rot}_v(K)$, also called Maslov index, can be defined
just as in the case where $K$ is Legendrian in a contact structure. Now:

\begin{lem}\label{wind:rot}
If $M$ is a homology sphere,
$K_0$ and $K_1$ are pseudo-Legendrian in $(M,v)$ and framed isotopic, then
$w(K_1,K_0)=\frac{1}{2}({\rm rot}_v(K_1)-{\rm rot}_v(K_0))$.
\end{lem}
(Concerning the statement, note that ${\rm rot}_v(K_1)-{\rm rot}_v(K_0)$ 
must be even if $K_0$ and $K_1$ are framed-isotopic, otherwise
one of $\{K_0,K_1\}$ would lift to a closed path in a spin structure on $M$, and
the other one would not: a contradiction. A proof is easily
obtained by isotoping $K_1$ as stated in Proposition~\ref{trace:like}.)

Lemma~\ref{wind:rot} gives another proof of the fact that $w(K_1,K_0)\in\mz$
can be defined when $M$ is a homology sphere. Moreover, it could be used 
to show goodness of knots in a homology sphere by a more direct argument than
that given in~\cite{first:paper}. We conclude this paragraph by showing
the result stated in the introduction and asking a question which 
naturally arises from it.

\dim{appl:wind:prop}
Equivalence of (\ref{iso:intro:point}) and (\ref{zero:wind:intro:point})
comes from Corollary~\ref{zero:wind}. Equivalence of (\ref{zero:wind:intro:point})
and (\ref{same:maslov:intro:point}) comes from Lemma~\ref{wind:rot}.
Equivalence of (\ref{same:maslov:intro:point}) and (\ref{same:torsion:intro:point})
follows from Corollary\ref{zero:wind} and the fact that the first homology group of
the complement of a knot in $M$ is infinite cyclic and generated by a meridian.
Equivalence of (\ref{same:maslov:intro:point}) and 
(\ref{imm:intro:point}) is an application of Gromov's $h$-principle
(see~\cite{vlad}).\finedim{appl:wind:prop}

\begin{que} Let $(M,v)$ be an arbitrary combed manifold, 
let $K_0$ and $K_1$ be pseudo-Legendrian in $(M,v)$ and framed-isotopic, and assume
that they are homotopic through pseudo-Legendrian immersions. Does this imply that
$w(K_1,K_0)\cdot[\mu_{K_0}]=0$? (We do not think that the opposite implication
can be true in general, in particular when $[\mu_{K_0}]$ has finite order.)
\end{que}

\paragraph{Absolute winding number}
We concentrate in this paragraph on fields $v$ such that $\ee(v^\perp)=0$, 
where $\ee$ denotes the Euler class and the choice of the metric is of course immaterial.
Condition $\ee(v^\perp)=0$ is equivalent to the existence of another non-vanishing field $x$
always transversal to $v$. Since the ambient manifold is oriented, this is also equivalent
to the fact that $v$ extends to a framing $(v,x,y)$, {\em i.e.} a global trivialization
of the tangent bundle to $M$. Assume now that $K$ is an oriented knot transversal to $v$.
Then, taking the projection of the tangent vector to $K$ on the
unit sphere of the $(x,y)$-plane, and computing the degree, we can define a rotation number
${\rm rot}_{(v,x)}(K)$.

\begin{rem} {\em ${\rm rot}_{(v,x)}(K)$ is invariant under simultaneous homotopy $(v_t,x_t)$ and
isotopy $(K_t)$ such that $x_t$ and $v_t$ are transversal to $v_t$ for all $t$. Moreover
${\rm rot}_{(v,x)}(K)$ is independent of $x$ when $M$ is a homology sphere,
and it equals the Maslov index already discussed above.}\end{rem}

Assume now that $K_0$ and $K_1$ are both transversal to $v$. Within the proof of
Proposition~\ref{trace:like} we have shown that $K_1$ can be isotoped through knots
transversal to $v$ to a knots which differs from $K_0$ by double curls only.

\begin{rem} {\em ${\rm rot}_{(v,x)}(K_1)-{\rm rot}_{(v,x)}(K_0)$ is independent of $x$
and equals twice the number of double curls by which $K_0$ and $K_1$ differ, up
to isotopy transversal to $v$.}\end{rem}

The previous remark shows that the relative winding number is well-defined (without
any assumption on the knots) if one
restricts to knots transversal to a given $v$ with $\ee(v^\perp)=0$, and one views
the knots up to {\em isotopy transversal to $v$} (as opposed to pseudo-Legendrian isotopy,
which involves also a homotopy of $v$). More on the difference between transversal isotopy
and pesudo-Legendrian isotopy will be said below.

\begin{rem}
{\em Combining the previous two remarks one gets yet another proof 
that the relative winding number is well-defined
in $\mz$ up to pseudo-Legendrian isotopy in a homology sphere.}\end{rem}

\paragraph{Finite-order invariants}
We formally state and motivate in this paragraph the conjecture announced in the
introduction. Let $\xi$ be an oriented contact structure on $M$ (which we assume
to be closed by simplicity), and let $v$ be a 
field positively transversal to $\xi$. 
Consider the spaces
$\Leg(M,\xi)$, $\Pleg^{{\rm weak}}(M,v)$
and $\Fram(M)$ of $\xi$-Legendrian, $v$-transverse, and framed knots in $M$, with the
appropriate equivalence relations (namely $\xi$-Legendrian, pseudo-Legendrian, and
framed isotopy).
Enlarge these spaces by allowing immersions of $S^1$ rather than
embeddings, and take path-connected components $\lgot$, $\pgot$ and $\fgot$,
with $\lgot\subset\pgot\subset\fgot$. (Concerning $\Pleg$, note that a path is a family
$(K_t,v_t)_{t\in[0,1]}$ with $v_0=v_1=v$.)
Given an Abelian group $A$ one can define, using
the customary Vassiliev-Goussarov skein relations, the spaces
$V_\lgotpicc^n(A)$, $V_\pgotpicc^n(A)$ and $V_\fgotpicc^n(A)$ of $A$-valued
order-$n$ invariants under Legendrian, pseudo-Legendrian and framed isotopy respectively.
Since a Legendrian isotopy is pseudo-Legendrian, and a pseudo-Legendrian isotopy is framed,
using restrictions we get a commutative diagram of homomorphisms:
$$\matrix{
V_\fgotpicc^n(A) & 
\mathop{\longrightarrow}\limits^{\phi^n_{\fgotpicc,\pgotpicc}} & 
V_\pgotpicc^n(A) \cr
& & \cr
& 
\scriptstyle{\phi^n_{\fgotpicc,\lgotpicc}}\searrow\ \ \ \ &
\ \ \downarrow\scriptstyle{\phi^n_{\pgotpicc,\lgotpicc}} \cr
& & \cr
&
& 
V_\lgotpicc^n(A).\cr}$$	
Tchernov's arguments~\cite{vlad} imply that
all three $\phi$'s are always injective, 
and our conjecture is that $\phi^n_{\pgotpicc,\lgotpicc}$ is always an isomorphism. 
By~\cite{vlad} again,
the conjecture is equivalent to showing that every finite-order Legendrian invariant
is automatically invariant also under pseudo-Legendrian isotopy.
The generalized Fuchs-Tabachnikov theorem (see~\cite{futa} and~\cite{vlad})
states that $\phi^n_{\fgotpicc,\lgotpicc}$ is an isomorphism in many cases
({\em e.g.}~if $M$ is a homology sphere), so $\phi^n_{\pgotpicc,\lgotpicc}$
is also an isomorphism in these cases. Tchernov has provided the only known examples
in which $\phi^n_{\fgotpicc,\lgotpicc}$ is not an isomorphisms, namely he has exhibited
elements of $V^n_\lgotpicc(A)$ which do not lift to $V^n_\fgotpicc(A)$.
Our impression is that these elements do lift to $V^n_\pgotpicc(A)$, which would
imply that $\phi^n_{\pgotpicc,\lgotpicc}$ is indeed an isomorphism in all known cases.
Truthness of our conjecture would imply that Legendrian \foi s are only sensitive to the
homotopy class of a contact structure, and in particular that they cannot capture
tightness.

\subsection{Pseudo-Legendrian vs.~Legendrian knots}\label{pseudo:section}
After the work of Eliashberg~\cite{eliash}, we know that on a closed manifold
an overtwisted contact structure is determined up to isotopy by its homotopy
class as a plane field. We discuss in this section the extent to which
this fact extends in presence of a pseudo-Legendrian link.
We start with an open question which arises from the results of 
Section~\ref{leg:calc} and will lead us to the connection with
overtwisted contact structures.

\paragraph{Fixed vs.~variable spine for pseudo-Legendrian links}
We will adopt in this paragraph the viewpoint which allows to
dismiss automorphisms of manifolds, fixing $M$ and considering spines
and moves embedded in $M$, as explained after the statement of
Proposition~\ref{rigid:spines}.

Theorems~\ref{comb:calc:teo} and~\ref{leg:calc:teo} 
and Proposition~\ref{fixed:P:prop} leave the following question open:
given an embedded spine $P\subset M$ representing 
a concave combing on $M$, what intrinsic topological
object is represented by $\cont^1$-diagrams on $P$ up to pseudo-Legendrian Turaev moves? 
Let us introduce some notation to formalize the situation. We denote by
$\dd^{\Pleg}(P)$ the set of equivalence classes of $\cont^1$-diagrams on $P$ up 
to pseudo-Legendrian Turaev moves. We also fix a representative 
$v$ of the combing carried by $P$ (so, $v$ is positively transversal to $P$
and restricts to $\bthtr$ on the complement of $P$).
We consider now the set of links in
$M$ transversal to $v$, and we denote by $\Pleg^{\rm weak}(M,v)$
the quotient space under the relation of existence of a pseudo-Legendrian isotopy,
{\em i.e.}~a path $(L_t,v_t)_{t\in[0,1]}$
as usual, with $v_0=v_1=v$. We also denote by 
$\Pleg^{\rm strong}(M,v)$ the (bigger) quotient obtained by forcing
$(v_t)$ to be constant. So $\Pleg^{\rm strong}(M,v)$ is just the set of
equivalence classes of links transversal to $v$. Using 
Proposition~\ref{fixed:P:prop} one sees
that the operation of turning a diagram into a link defines a bijection
$$\psi^{\rm strong}:\dd^\Pleg(P)\to\Pleg^{\rm strong}(M,v).$$
(This is not quite the content of Proposition~\ref{fixed:P:prop},
because here $(M,v)$ is $(\hatM(P),\hatv(P))$ rather that $(M(P),v(P))$,
but a link isotopy can be modified to avoid a $\bthtr$, and the conclusion follows.)

Bijectivity of $\psi^{\rm strong}$
is significant if one imagines to have started with the pair $(M,v)$, and
to have constructed $P$ from a normal section of $v$, 
as in Proposition~\ref{section:then:spine}. 
It is however less significant if one assumes only $P$ to be given from the beginning,
because in this case $v$ is actually well-defined only up
to homotopy, and fixing a representative looks artificial. The natural map to consider
is in this case 
$$\psi^{\rm weak}:\dd^\Pleg(P)\to\Pleg^{\rm weak}(M,v),$$
obtained by composition with the projection $\Pleg^{\rm strong}(M,v)
\to\Pleg^{\rm weak}(M,v)$. This map is of course surjective, and
one can ask whether it is injective or not. Some remarks are in order:
\begin{enumerate}
\item Theorem~\ref{leg:calc:teo} implies that if $\psi^{\rm weak}(D)=
\psi^{\rm weak}(D')$ then there exists a circular sequence $P=P_0\to P_1\to\cdots\to P_n=P$
of sliding-MP-moves and diagrams $D_i,D'_i\in\dd(P_i)$ with $D_0=D$, $D'_n=D'$,
$D_i\to D'_i$ a pseudo-Legendrian Turaev move, and $D_{i+1}$ the companion of $D'_i$ through
$P_i\to P_{i+1}$. Checking the injectivity of $\psi^{\rm weak}$ corresponds to 
the (purely combinatorial) question whether such a 
sequence $(P_i,D_i)$ can be replaced by one with constant $P_i$.
\item Using Theorem~\ref{comb:calc:teo} and the fact that $\cont^1$ diagrams can be followed
through sliding-MP-moves, one sees quite easily that injectivity of 
$\psi^{\rm weak}$ actually depends only on the combing carried by $P$, not on
$P$ itself.
\item Injectivity of $\psi^{\rm weak}$ is equivalent to injectivity of
the projection $\Pleg^{\rm strong}(M,v)
\to\Pleg^{\rm weak}(M,v)$, a purely topological question.
Injectivity of this projection may appear very unlikely at first sight,
since it basically corresponds to the fact that a homotopy can be replaced by an
isotopy. However one can remark that injectivity of projection depends
only on the homotopy class of $v$, rather than $v$ itself, so one
can assume that $v$ is transversal to an overtwisted contact structure.
For overtwisted structures, after the work of Eliashberg~\cite{eliash}, 
it is indeed true that homotopy implies isotopy, but the presence of
the link of course somewhat modifies the situation. 
We will expound this theme in the next paragraph.
\end{enumerate}

\paragraph{Overtwisted structures and overtwisted knot complements}
We fix in this paragraph an overtwisted contact structure $\xi$ on $M$ (which we
assume to be closed by simplicity) and a field $v$ positively transversal to $\xi$.
We will denote by $\Leg(M,\xi)$ the space of Legendrian links in $(M,\xi)$ 
up to Legendrian isotopy. In~\cite{first:paper}, having also in mind the
facts mentioned in the previous paragraph, we put forward the question of whether
the natural map 
$$\Leg(M,\xi)\to\Pleg^{\rm weak}(M,v)$$
is a bijection. A fact implying that this map is not injective in some cases
was recently communicated to us by E.~Giroux~\cite{giroux:mail}. 
He was able to construct triples
$(M,\xi,K)$ where $\xi$ is overtwisted, $K$ is $\xi$-Legendrian, and 
$\xi\ristr{M\setminus K}$ is tight. Let us apply a Lutz twist away from $K$
to get a new structure $\xi'$ such that $\xi'$ is 
homotopic to $\xi$ as a plane field on $M$, and $\xi'\ristr{M\setminus K}$ is overtwisted.
Using Eliashberg's classification~\cite{eliash} we consider $\varphi\in{\rm Diff}_0(M)$
such that $\xi'=\varphi^*(\xi)$, and define $K'=\varphi(K)$. By construction
$K$ and $K'$ have the same image in $\Pleg^{\rm weak}(M,v)$, but of course
they are inequivalent in $\Leg(M,\xi)$.

To avoid the phenomenon discovered by Giroux we consider in
$\Leg(M,\xi)$ the subset $\Leg^{\rm OT}(M,\xi)$ given by links whose complement
is overtwisted. We start by showing:

\begin{prop}\label{OT:surj}
The natural map $\Leg^{\rm OT}(M,\xi)\to\Pleg^{\rm weak}(M,v)$ is surjective.
\end{prop}

\dim{OT:surj}
Let $L$ be transversal to $v$, and fix a metric on $M$. 
Let $\eta$ be a positive contact structure near $L$ with $\eta=v^\perp$ on $K$
(such an $\eta$ is unique up to isomorphism). Extend $\eta$ to any plane field
homotopic to $v^\perp$ (and hence to $\xi$) on $M$. So $\eta$ is a plane distribution
which has a contact zone, and $L$ lies in this contact zone. The technique of 
Eliashberg~\cite{eliash} now allows to homotope $\eta$ away from its contact zone
to an overtwisted contact structure $\xi'$. 
The resulting $\xi'$ is now isotopic to $\xi$, again by Eliashberg's result. 
If $\varphi\in{\rm Diff}_0(M)$ and 
$\xi'=\varphi^*(\xi)$ we define $L'=\varphi(L)$. By construction 
$(L',v)$ is pseudo-Legendrian isotopic to $(L,v)$, and surjectivity is proved.
\finedim{OT:surj}

We cannot presently state whether the map
$\Leg^{\rm OT}(M,\xi)\to\Pleg^{\rm weak}(M,v)$ is injective or not in general. We only
give a partial argument (based on the techniques of Eliashberg again),
and mention where the difficulty arises.
Assume that $L_0$ and $L_1$ are $\xi$-Legendrian with
overtwisted complements and define equivalent pseudo-Legendrian links.
Then there exists a continuous family $(L_t,\xi_t)_{t\in[0,1]}$, where $\xi_0=\xi_1=\xi$
but $\xi_t$ is only a plane field for $t\neq 0,1$.
Eliashberg's contactization methods for homotopies,
together with the uniqueness of contact structures in the neighbourhood of
Legendrian links, should in our opinion allow to 
replace such a $(\xi_t)_{t\in[0,1]}$ by another one in which each $\xi_t$ 
is a contact structure (and still contains $L_t$ as a Legendrian link).
Applying Gray's theorem we get an isotopy $(\varphi_t)_{t\in[0,1]}$
such that $\xi_t=\varphi^*_t(\xi_0)$. Setting $\tilde L_t=\varphi_t(L_t)$
we get a Legendrian isotopy between $L_0$ and $\varphi_1(L_1)$. 
The question whether $\varphi_1(L_1)$ is automatically Legendrian isotopic to $L_1$,
at least for some classes of manifolds, now depends on the analysis of the
group ${\rm Aut}(M,\xi)\cap{\rm Diff}_0(M)$, 
which we leave unsettled for the time being.

\vspace{1cm}

\hspace{8cm} benedett@dm.unipi.it

\hspace{8cm} petronio@dm.unipi.it

\hspace{8cm} Dipartimento di Matematica

\hspace{8cm} Via F.~Buonarroti, 2

\hspace{8cm} I-56127, PISA (Italy)


\begin{thebibliography}{BePe3}

\bibitem[1]{manuscripta} {\sc R.~Benedetti, C.~Petronio}, {\it A finite
graphic calculus  for $3$-manifolds}, Manuscripta Math. {\bf 88} (1995),
291-310.

\bibitem[2]{lnm}  {\sc R.~Benedetti, C.~Petronio}, ``Branched Standard
Spines of 3-Manifolds'', Lecture Notes in Math. n. 1653, Springer-Verlag,
Berlin-Heidelberg-New York, 1997.

\bibitem[3]{contspin} {\sc R.~Benedetti, C.~Petronio},  {\it Branched spines
and contact structures on $3$-manifolds}, Ann. Mat. Pura Appl. (to appear).

\bibitem[4]{first:paper} {\sc R.~Benedetti, C.~Petronio},  
{\it Reidemeister torsion of 3-dimensional Euler structures
with simple boundary and Legendrian links}, preprint {\tt Math.GT/9907184}, submitted.

\bibitem[5]{benvenuti} {\sc S. Benvenuti}, {\it Hopf algebras and invariants of combed and 
framed 3-manifolds}, In: ``Knots in
Hellas '98,'' to appear on a special issue of J. Knot Theory Ramif.

\bibitem[6]{casler} {\sc B.~G.~Casler}, {\it An imbedding theorem for
connected $3$-manifolds with boundary}, Proc. Amer. Math. Soc. {\bf 16} (1965),
559-566.

\bibitem[7]{eliash} {\sc Ya.~Eliashberg}, {\it Classification of 
overtwisted contact structures}, Invent. Math. {\bf 98} (1989), 623-637.

\bibitem[8]{futa} {\sc D.~Fuchs, S.~Tabachnikov}, 
{\it Invariants of Legendrian and transverse knots in the standard contact space},
Topology {\bf 36} (1997), 1025-1053.

\bibitem[9]{giroux:mail} {\sc E.~Giroux}, private communication,
December 1999.

\bibitem[10]{ishii} {\sc I.~Ishii}, {\it Moves for flow-spines and
topological invariants of $3$-manifolds}, Tokyo J. Math. {\bf 15} (1992),
297-312.

\bibitem[11]{matv:mossa} {\sc S.~V.~Matveev}, {\it Transformations of special
spines and the Zeeman conjecture}, Math. USSR-Izv. {\bf 31} (1988), 423-434.

\bibitem[12]{meta}{\sc G.~Meng, C.~Taubes},
{\it $\underline{SW}\,=\,$Milnor torsion}, Math. Res. Lett. {\bf 3} (1996), 661-674

\bibitem[13]{piergallini} {\sc R.~Piergallini}, {\it Standard moves for  
standard polyhedra and spines}, Rendiconti Circ. Mat. Palermo {\bf 37}, suppl.
18 (1988), 391-414.

\bibitem[14]{vlad} {\sc V.~Tchernov}, {\it Finite order invariants of Legendrian,
transverse, and framed knots in contact 3-manifolds}, preprint 1999.

\bibitem[15]{trace} {\sc B.~Trace}, {\it On the Reidemeister
moves of a classical knot}, Proc. Amer. Math. Soc. {\bf 89} (1983), 722-724.

\bibitem[16]{turaev:ombre} {\sc V.~G.~Turaev}, ``Quantum Invariants 
of Knots and 3-Manifolds'', de Gruyter stud. in Math. 18, Berlin-New York, 1994.

\bibitem[17]{turaev:Euler} {\sc V.~G.~Turaev}, {\it Euler structures, nonsingular
vector fields, and torsion of Reidemeister type}, Math. USSR-Izv. {\bf 34}
(1990), 627-662.

\bibitem[18]{turaev:spinc} {\sc V.~G.~Turaev}, {\it Torsion invariants of
Spin$^{\it c}$-structures on $3$-manifolds}, Math. Res. Lett. {\bf 4} (1997), 679-695.

\bibitem[19]{turaev:nuovo} {\sc V.~G.~Turaev}, {\it A combinatorial
formulation for Seiberg-Witten invariants of $3$-manifolds},
Math. Res. Lett. {\bf 5} (1998), 583-598.


\end{thebibliography}
\end{document}